\documentclass[11pt]{article}

\usepackage[margin=1.2in]{geometry}

\usepackage{amsmath,amssymb,amsfonts,graphicx,amsthm,nicefrac,mathtools,bm}

\usepackage{hyperref}
\usepackage{natbib}
\usepackage{bbm}
\usepackage{color}
\usepackage{xspace}
\newtheorem{theorem}{Theorem}

\newtheorem{lemma}{Lemma}
\newtheorem{proposition}{Proposition}

\newtheorem{assumption}{Assumption}
\theoremstyle{definition}

\theoremstyle{remark}

\makeatletter
\newtheorem*{rep@theorem}{\rep@title}
\newcommand{\newreptheorem}[2]{%
\newenvironment{rep#1}[1]{%
\def\rep@title{#2 \ref{##1}}%
\begin{rep@theorem}}%
{\end{rep@theorem}}}
\makeatother
\newreptheorem{theorem}{Theorem}
\newreptheorem{lemma}{Lemma}
\newreptheorem{corollary}{Corollary}
\newreptheorem{proposition}{Proposition}

\DeclareMathOperator{\sign}{sign}

\DeclareMathOperator*{\argmin}{arg\,min}

\newcommand{\One}[1]{{\mathbbm{1}}\left\{{#1}\right\}}
\newcommand{\inner}[2]{\langle{#1},{#2}\rangle} % Inner product
\newcommand{\norm}[1]{\lVert{#1}\rVert}

\newcommand{\PP}[1]{\mathbb{P}\left\{{#1}\right\}} % Probability
\newcommand{\EE}[1]{\mathbb{E}\left[{#1}\right]} % Expectation

\newcommand{\EEst}[2]{\mathbb{E}\left[{#1}\ \middle| \ {#2}\right]} % Conditional expectation
\renewcommand{\O}[1]{\mathcal{O}\left({#1}\right)}
\def\R{\mathbb{R}}

\newcommand{\ident}{\mathbf{I}}

\newcommand\independent{\protect\mathpalette{\protect\independenT}{\perp}}
\def\independenT#1#2{\mathrel{\rlap{$#1#2$}\mkern2mu{#1#2}}}

\newcommand{\iidsim}{\stackrel{\mathrm{iid}}{\sim}}

\newcommand{\eps}{\varepsilon}

\newcommand{\xh}{\tilde{x}}
\newcommand{\yh}{\tilde{y}}
\newcommand{\uh}{\tilde{u}}
\newcommand{\biginner}[2]{\left\langle{#1},{#2}\right\rangle}
\newcommand{\dom}{\textnormal{dom}}
\newcommand{\Lcal}{\mathcal{L}}

\newcommand{\loss}{\textnormal{Loss}}
\DeclareMathOperator{\qexp}{qexp}

\newcommand{\normal}{\mathcal{N}}
\DeclareMathOperator{\vect}{vec}

\usepackage{algorithm,algorithmicx,algpseudocode}
\algdef{SE}[FORUNTIL]{for}{until}[1]{\algorithmicfor \ #1 \algorithmicdo}[1]{\algorithmicuntil\ #1}%

\title{Convergence for nonconvex ADMM,\\ with applications to CT imaging}
\author{Rina Foygel Barber\thanks{\texttt{rina@uchicago.edu} ; Department of Statistics, University of Chicago} \ and Emil Y.~Sidky\thanks{\texttt{sidky@uchicago.edu} ; Department of Radiology, University of Chicago}}
\date{\today}

\begin{document}
\maketitle

\begin{abstract}
The alternating direction method of multipliers (ADMM) algorithm is a powerful and flexible tool for complex optimization problems
of the form $\min\{f(x)+g(y) : Ax+By=c\}$. ADMM exhibits robust empirical performance
across a range of challenging settings including nonsmoothness and nonconvexity of the objective functions $f$ and $g$, and
 provides a simple and natural approach to the inverse
problem of image reconstruction for computed tomography (CT) imaging.
 From the theoretical point of view,
existing results for convergence in the nonconvex setting generally assume smoothness in at least one of the component functions in the objective.
In this work, our new theoretical results provide
convergence guarantees under a restricted strong convexity assumption without requiring smoothness or differentiability, while still
allowing differentiable terms to be treated approximately if needed. We validate these theoretical results empirically, with a simulated example
where both $f$ and $g$ are nondifferentiable---and thus outside the scope of existing theory---as well as a simulated CT image reconstruction
problem.
\end{abstract}

\section{Introduction}
In this work, we consider optimization problems of the form
\begin{equation}\label{eqn:optim_intro}\textnormal{Minimize }f(x) + g(y) \textnormal{ subject to the constraint that } Ax + By = c.\end{equation}
Problems of this form arise in many applications throughout the physical and biological sciences. In particular, we are interested in optimization problems
pertaining to computed tomography (CT) imaging, which, as we will see later on, can often be expressed in this type of formulation.

Solving the optimization problem~\eqref{eqn:optim_intro} can be computationally challenging even when the functions $f$ and $g$ are both
convex. Challenges in the convex setting may include high dimensionality of the variables $x$ and $y$, nondifferentiability of $f$ and/or $g$,
or poor conditioning of the linear transformations $A,B$ or the functions $f,g$. 
If one or both functions are nonconvex, this brings an additional level of difficulty to the optimization problem. 

In this work, we study a linearized form of the 
alternating directions method of multipliers (ADMM) algorithm, in the setting 
where $f$ and $g$ may both be nonconvex and nonsmooth. While variants of this algorithm are very well known in the literature,
existing theoretical results have typically been restricted to narrower settings (e.g., assuming that at least one of the two functions $f$, $g$ must
be smooth), and thus cannot be applied to guarantee convergence for many settings arising in modern high dimensional 
optimization and data analysis. 

\paragraph{Outline} In Section~\ref{sec:setting}, we describe the method of nonconvex ADMM with linear approximations, and review known results
in the literature on the convergence properties of this type of algorithm in various settings.
In Section~\ref{sec:theory} we present our new convergence result, which addresses a more flexible setting allowing both $f$ and $g$
to be potentially nonconvex and nonsmooth.
We demonstrate the performance of the algorithm on a simple simulated quantile regression estimation problem in Section~\ref{sec:examples},
and present an application to computed tomography (CT) imaging in Section~\ref{sec:CT}. 
Finally, some future directions and implications of this work are discussed
in Section~\ref{sec:discussion}. Some proofs and additional technical details are deferred to the Appendix.

\section{Setting and background}\label{sec:setting}
Consider the optimization problem
\begin{equation}\label{eqn:optim}\textnormal{Minimize }f(x) + g(y) \ : \ x\in\R^d,y\in\R^m\textnormal{ such that } Ax + By = c\end{equation}
where the functions $f$ on $\R^d$ and $g$ on $\R^m$ are potentially nonconvex and/or nondifferentiable, while $A\in\R^{k\times d}$, $B\in\R^{k\times m}$, and $c\in\R^k$
define linear constraints on the variables.
In this work, we will consider functions $f$ and $g$ that can be decomposed as
\[f(x) = f_c(x) + f_d(x), \quad g(y) = g_c(y) + g_d(y)\]
where $f_c$ is convex (possibly nondifferentiable) and $f_d$ is twice differentiable (possibly nonconvex), and similarly for $g_c$ and $g_d$.
This decomposition allows us to take linear approximations to the differentiable terms $f_d$ and $g_d$, where needed, 
to ensure simple calculations for each update step of our iterative algorithm.

We will assume that $f_c$ and $g_c$ are {\em proper functions}. Formally, this means that we can write 
\[f_c: \R^d \rightarrow \R\cup \{+\infty\},\]
with nonempty domain $\dom(f_c) := \{x\in\R^d: f(x)<+\infty\}$ (and similarly for $g_c$). We also assume that $f_c$ and $g_c$ are lower semi-continuous.
The differentiable component $f_d$ is assumed to be defined on all of $\R^d$, i.e.,
\[f_d:\R^d\rightarrow\R,\]
and similarly for $g_d$ on $\R^m$. Putting these assumptions together, we see that $f$ and $g$ are also proper functions,
with domains $\dom(f)=\dom(f_c)$ and $\dom(g)=\dom(g_c)$ (note that convexity of $f_c,g_c$ ensures that these domains are also convex). 
Finally, we assume that the feasible set
\[\big(\dom(f)\times \dom(g)\big) \cap \left\{(x,y)\in\R^d\times\R^m : Ax+By=c\right\}\]
is nonempty. We will say that a point $(x,y)$ is {\em feasible} for this optimization problem if it lies in 
this feasible set, i.e., $x\in\dom(f)$, $y\in\dom(g)$, and the constraint
$Ax+By=c$ is satisfied. 

\subsection{Background and prior work}
\subsubsection{ADMM for convex optimization problems}
The alternating directions method of multipliers (ADMM) algorithm is a method
for solving problems of the form~\eqref{eqn:optim}. 
It was developed initially for the setting where $f$ and $g$ are both convex,
and operates by reformulating the optimization problem~\eqref{eqn:optim}
with an augmented Lagrangian,
\[\min_{x,y}\max_u\left\{\Lcal_{\Sigma}(x,y,u)\right\},\]
where the augmented Lagrangian is defined as
\begin{equation}\label{eqn:augLagr}\Lcal_{\Sigma}(x,y,u) = f(x) + g(y) + \inner{u}{Ax+By-c} + \frac{1}{2}\norm{Ax+By-c}^2_\Sigma,\end{equation}
for some positive definite penalty matrix $\Sigma\succ 0$. (Most commonly, $\Sigma$ is taken
to be a multiple of the identity.) See \citet{boyd2011distributed}
for a review of the motivation and performance of ADMM for the convex setting, including the long history of this algorithm and
many of its variants.

The ADMM algorithm solves this optimization problem as follows: initializing at some $x_0,u_0,y_0$, for all $t\geq 0$ we run the steps:
\begin{equation}\label{eqn:ADMM_iters_basic}\begin{cases}
x_{t+1} = \argmin_x\left\{\Lcal_{\Sigma}(x,y_t,u_t)\right\},\\
y_{t+1} = \argmin_y\left\{\Lcal_{\Sigma}(x_{t+1},y,u_t)\right\},\\
u_{t+1} = u_t + \Sigma(Ax_{t+1}+By_{t+1} - c).
\end{cases}\end{equation}

\paragraph{Adding step size matrices}
In some cases, adding step size matrices $H_f\succeq 0$ for the $x$ update and $H_g\succeq 0$ for the $y$ update
can improve the convergence behavior and/or may allow for easier calculation of the update steps:
\begin{equation}\label{eqn:ADMM_iters_stepsize}\begin{cases}
x_{t+1} = \argmin_x\left\{\Lcal_{\Sigma}(x,y_t,u_t) + \frac{1}{2}\norm{x-x_t}^2_{H_f}\right\},\\
y_{t+1} = \argmin_y\left\{\Lcal_{\Sigma}(x_{t+1},y,u_t) + \frac{1}{2}\norm{y-y_t}^2_{H_g}\right\},\\
u_{t+1} = u_t + \Sigma(Ax_{t+1}+By_{t+1} - c).
\end{cases}\end{equation}
(Here $\succeq$ denotes the positive semidefinite ordering on matrices, i.e., $H_f\succeq 0$ means that $H_f$ 
is positive semidefinite.)

In many cases, choosing $H_f$ so that $D_f:=H_f + A^\top \Sigma A$ is diagonal, or is a multiple of the identity, may be convenient for calculating the $x$ update step---this is 
because the $x$ update step is a minimization problem of the form $\argmin_x\left\{f(x) + \frac{1}{2}x^\top D_f x - x^\top v_t\right\}$, where $v_t$
is a vector that depends on the previous iteration.
Specifically, this type of choice for $H_f$ can be helpful when the function $f$ separates over the entries of $x$, $f(x) = \sum_i f_i(x_i)$,
so that now the $x$ update step separates completely over the entries of $x$. 
Another setting
where this type of modification is commonly used is when $f$ is equipped with an inexpensive proximal map
 (the map $z\mapsto \argmin\{f(x)+\frac{1}{2}\norm{x-z}^2_2\}$)---for example, the $\ell_1$ norm, $f(x) = \norm{x}_1$, or the (squared) $\ell_2$ norm, $f(x)=\norm{x}^2_2$,
are both commonly used regularization functions that have simple proximal maps.
(Without the matrix $H_f$, the $x$ update step is of the form  $\argmin_x\left\{f(x) + \frac{1}{2}x^\top A^\top\Sigma A x - x^\top v_t\right\}$, which may be substantially
more challenging to compute if $A^\top\Sigma A$ is a dense matrix.)
Similarly we may choose $H_g$ with these types of considerations in mind for the $y$ update step.  For further
details, see \citet[Eqn.~(17)]{wang2014bregman}, where this type of modification is referred to as a ``linearization''
of the quadratic penalty term.

This type of modification of ADMM is closely linked to related algorithms
for composite optimization problems of the form $f(Ax) + g(x)$, studied via primal-dual methods by, e.g.,  \citet{chen1994proximal,chambolle2011first,he2012convergence,valkonen2014primal},
among many others, and has been applied to convex versions of the CT image reconstruction problem (see, e.g., \citet{nien2014fast}).

\paragraph{Linear approximations} For many optimization problems, even with the modification of a step size matrix as in~\eqref{eqn:ADMM_iters_stepsize} above,
it may still be challenging to compute the $x$ update step if the function $f$ is difficult to minimize (and similarly, the $y$ step with the function $g$).
 In particular, if the $x$ update step itself can only be solved with an iterative procedure, this type of ``inner loop'' will drastically slow down the convergence of ADMM.
 
An alternative is to replace the function $f$ with an approximation at each step. In particular, consider our earlier decomposition, $f=f_c+f_d$, where $f_c$ is convex
while $f_d$ is twice differentiable. Taking a linear approximation to $f_d$,
at the current iteration $x_t$, we can  approximate the function $f$ as
\[f(x)\approx f_c(x) + \big(f_d(x_t) + \inner{\nabla f_d(x_t)}{x - x_t}\big).\]
Although this inexact calculation of the $x$ update may lead to slower convergence in terms of the total number of iterations, this may be outweighed if
this approximation allows the cost of each single iteration to be substantially reduced.
We can make the analogous modification for the $y$ update step.
This  type of modification has been commonly used in both the convex  and nonconvex settings, particularly
in settings where $f$ itself is twice differentiable so we can take $f_d=f$ and $f_c\equiv 0$. 
For instance, \citet[Eqn.~(21)]{wang2014bregman} study this modification for the convex setting, where this type
of approach is referred to as ``linearization'' of the target function; see also the references
described below
for the nonconvex setting.

\begin{algorithm}[t]
\caption{ADMM with linear approximations}
\label{algo:main}
\begin{algorithmic}
\State \textbf{Input:} Functions $f = f_c+f_d$ and $g=g_c+g_d$, with $f_c,g_c$ convex, $f_d,g_d$ twice differentiable;\\
\quad\quad \quad\quad  matrices $A,B$; vector $c$; penalty matrix $\Sigma\succ 0$; step size matrices $H_f,H_g\succeq 0$.
\State \textbf{Initialize:} $x_0,y_0,u_0$.
\for{$t=0,1,2,\dots$} 
\begin{align*}&\textnormal{Update $x$:\quad}
x_{t+1} = \argmin_x \bigg\{f_c(x) + \inner{x}{\nabla f_d(x_t) + A^\top u_t} \\&\hspace{2in}{}+\frac{1}{2}\norm{Ax +By_t -c}^2_\Sigma + \frac{1}{2}\norm{x - x_t}^2_{H_f} \bigg\}.\\
&\textnormal{Update $y$:\quad}
y_{t+1} = \argmin_y\bigg\{g_c(y) +\inner{y}{\nabla g_d(y_t) + B^\top u_t}\\&\hspace{2in}{} +\frac{1}{2}\norm{Ax_{t+1}+By-c}^2_\Sigma + \frac{1}{2}\norm{y - y_t}^2_{H_g}\bigg\}.\\
&\textnormal{Update $u$:\quad}
u_{t+1} = u_t + \Sigma(Ax_{t+1}+By_{t+1} - c).
\end{align*}
\until{some convergence criterion is reached.}
\end{algorithmic}
\end{algorithm}
For completeness, Algorithm~\ref{algo:main}
presents this modified form of ADMM (combining both linear approximations to $f_d$ and $g_d$, and the addition of step size matrices described
above). This is the version of the algorithm that we will study in our work.

\subsubsection{Nonconvex ADMM}\label{sec:prior_work_nonconvex}
Next we turn to the nonconvex setting, where the functions $f$ and/or $g$ are no longer required to be convex.
In many optimization problems, the ADMM algorithm (possibly with the addition of step size matrices $H_f,H_g$ and/or linear approximations to $f_d,g_d$)
has been observed to perform well, converging successfully and avoiding issues such as saddle points or local minima.
The convergence properties in a nonconvex setting have been studied extensively. For example, 
\citet{wang2014convergence,magnusson2015convergence,hong2016convergence,guo2017convergence,wang2018convergence,wang2019global,themelis2020douglas} study the performance of ADMM 
with $f$ and $g$ update steps calculated exactly (in some cases, extending the algorithm to handle more than two variable blocks), while  \citet{li2015global,lanza2017nonconvex,jiang2019structured,liu2019linearized}
study the algorithm with linear approximations to (parts of) $f$ and/or $g$. All of these works prove results of one of the two following types:
\begin{itemize}
\item Assume that either $f$ or $g$ is differentiable and has a Lipschitz gradient, and establish convergence guarantees;
\item Assume that the algorithm converges (or, more weakly, assume only that the dual variable $u_t$ converges),
and establish optimality properties of the limit point.
\end{itemize}
It is important to note that neither type of existing result verifies that convergence is guaranteed in a nonconvex
setting where both $f$ and $g$
are nondifferentiable. 

A different type of nonconvexity that is studied in the literature is where $f$ and $g$ are both convex, but the constraint on $(x,y)$ is nonconvex (e.g., $y=A(x)$ for a nonlinear
operator $A$); this type of problem is studied by \citet{valkonen2014primal,ochs2015iteratively}, among others.
\citet{bolte2018nonconvex} allow for nonconvexity both in the functions ($f$ and/or $g$) and in the constraint on $(x,y)$; as
with many of the methods above, the results of this paper require that either $f$ or $g$ is differentiable and has a Lipschitz gradient.

\subsubsection{The MOCCA algorithm}\label{sec:prior_work_mocca}

Our own earlier work on this problem \citep{barber2016mocca} proposed the Mirrored Convex/Concave algorithm (MOCCA), 
which solves problems of the form~\eqref{eqn:optim_intro}. 
At a high level, the MOCCA algorithm can be viewed as a version of Algorithm~\ref{algo:main} with a key modification:
rather than taking a new linear approximation to $f_d$ and $g_d$ at each iteration $t$ (i.e., computing the gradients $\nabla f_d(x_t)$ and $\nabla g_d(y_t)$),
the MOCCA algorithm requires an ``inner loop'', where we cycle  $L_t$ many times
 through the variable update steps before
re-calculating the linear approximations to $f_d$ and $g_d$. 

In \citep{barber2016mocca}, two versions of the MOCCA algorithm are proposed:
\begin{itemize}
\item The ``stable'' version \citep[Algorithm 2]{barber2016mocca}, where at each iteration $t$ of the outer loop, we run $L_t \gg1 $ many iterations of the inner loop, and require $L_t\rightarrow\infty$.
\item The ``simple'' version \citep[Algorithm 1]{barber2016mocca}, with no inner loop (or equivalently, with $L_t = 1$ for each $t$).
\end{itemize}
The theoretical guarantee given in \citep{barber2016mocca} proves a convergence result for the ``stable'' version.
To our knowledge, this was a unique result
 in that it ensured
convergence without requiring either $f$ or $g$ to have a Lipschitz gradient (in comparison
to the literature on ADMM in the nonconvex setting as discussed above), requiring instead a restricted strong convexity type condition (see
Section~\ref{sec:rsc} below). 
However, the theoretical result has the drawback of requiring the inner loop, with $L_t\rightarrow \infty$.
This requirement contradicts the empirical performance of the algorithm: the empirical results in  \citep{barber2016mocca}
actually implemented the ``simple'' version of MOCCA, with no inner loop, and the algorithm typically showed convergence 
even though no theoretical justification was known.

The ADMM algorithm studied in the present work, Algorithm~\ref{algo:main}, is in fact essentially equivalent to the ``simple'' version of MOCCA (with a few
changes in the details; e.g., in MOCCA, the matrix $B$ was required to be the identity). The novelty of the present work, then,
is not in the algorithm itself, but rather in the fact that the theoretical guarantees established in this paper apply to the actual algorithm being 
run in practice (Algorithm~\ref{algo:main}, or equivalently, the ``simple'' version of MOCCA), rather than applying only to a 
more computationally inefficient version of this algorithm (the ``stable'' version of MOCCA, as in the theoretical results of \citep{barber2016mocca}).

\subsection{Preview of new results}
In the present work, we establish a convergence guarantee for Algorithm~\ref{algo:main} in the nonconvex setting, with no ``inner loop'' needed in the theory, substantially closing
the gap between the theoretical results and our empirical observations for this algorithm.
As byproducts of this new analysis, we uncover an additional interesting finding that better explains the dependence of performance
on step size parameters. Moreover, our new work allows for a more direct connection to CT imaging---we are able to apply our algorithm,
exactly as defined and with no modifications, to simulated CT image reconstruction problems, obtaining very clean results.
(For real CT data, issues of scanner calibration, non-random noise, etc., require a more careful application of the algorithm,
which we address in separate work, but we mention here that this algorithm has been very successful on real CT data, 
e.g.,~\citet{rizzo2022material,schmidt2023constrained,rizzo2023experimental}.)

\section{Convergence guarantee}\label{sec:theory}

We will prove a convergence result under an additional condition requiring approximate convexity of the problem.

If the optimization problem were strongly convex, we would expect that our optimization algorithm would converge
to the unique minimizer---which, under strong convexity, would be a point satisfying first-order optimality conditions.
 In more challenging settings, however,  strong convexity may not hold and we will need to relax
our goal for convergence.

In this work, we will consider a setting where there is a feasible point $(\xh,\yh)$ that is \emph{approximately} first-order optimal,
around which the optimization problem satisfies a \emph{relaxed} version of strong convexity.
Since these conditions will only be required to hold approximately, the point $(\xh,\yh)$ may in general be nonunique; feasible
 points $(\xh',\yh')$ sufficiently close to $(\xh,\yh)$ might also satisfy the conditions. This is not a contradiction, however,
 since our theoretical results will only guarantee convergence to within some neighborhood of $(\xh,\yh)$.

In the remainder of this section, we will define our assumptions more formally and will state the theoretical guarantee, but we first need to
 review the definition of the subdifferential
in this nonconvex setting.

\subsection{Subdifferentials of $f$ and $g$}
\begin{figure}[t]\centering
\includegraphics[height=0.3\textwidth]{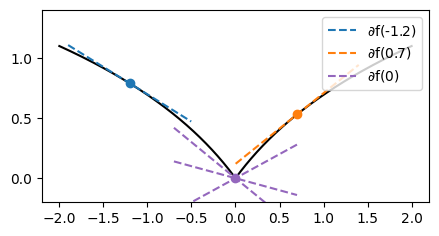}
\label{fig:subdifferential}\caption{Illustration of the subdifferential $\partial f(t)$, for the function $f(t) = \log(1+|t|)$.
For any $t\neq 0$, the function is differentiable at $t$, and the subdifferential is a singleton set containing only this derivative,
$\partial f(t) = \{f'(t)\} = \{ \textnormal{sign}(t)/(1+|t|)\}$. This is illustrated in the figure for two nonzero values of $t$. At $t=0$, the function is nondifferentiable,
and the subdifferential is given by $\partial f(0) = [-1,1]$. This is illustrated in the figure by showing several elements of $\partial f(0)$.}
\end{figure}

Since $f$ and $g$ are not necessarily convex, we pause here to define the notation $\partial f(x)$ and $\partial g(y)$,
which is a generalization of the usual subdifferential for convex functions. Here, for any $x\in\dom(f)$, we will use the definition
\[\partial f(x) = \left\{\xi : \lim_{t\rightarrow 0} \frac{f(x+tw)-f(x)}{t} \geq \inner{\xi}{w}\textnormal{ for all $w\in\R^d$}\right\}\]
and similarly for $g$. 
This definition is illustrated in Figure~\ref{fig:subdifferential}.

In particular, given the convex-plus-differentiable decomposition
$f=f_c+f_d$, we can write 
\[\partial f(x) = \{\xi + \nabla f_d(x) : \xi\in \partial f_c(x)\}\textnormal{\quad and \quad}\partial g(y) = \{\zeta+ \nabla g_d(y) : \zeta\in\partial g_c(y)\},\]
where $\partial f_c(x)$ and $\partial g_c(y)$ are the usual subdifferentials of the convex functions $f_c$ and $g_c$, i.e.,
for all $x\in\dom(f)$ we define
\[\partial f_c(x) = \left\{\xi : f(x+w)-f(x) \geq \inner{\xi}{w}\textnormal{ for all $w\in\R^d$}\right\},\]
and similarly for $g_c$. 

From this point on, for any $x\in\dom(f)$ and any $y\in\dom(g)$,
$\xi_x$ always denotes an element of $\partial f(x)$,  and $\zeta_y$ always denotes an element of $\partial g(y)$.

\subsection{Restricted strong convexity}\label{sec:rsc}

We will assume a restricted strong convexity (RSC) condition, which at a high level is a relaxation of imposing a strong convexity condition on the constrained optimization problem.
This type of convexity condition has been extensively studied in the high-dimensional statistics literature. For background,
the condition was proposed initially
by \citet{negahban2012unified}, and was studied by \citet{loh2015regularized} in the setting of nonconvex loss functions.
This type of condition is known to characterize many settings where accurate signal recovery is possible in spite of the ``curse of dimensionality'',
and over recent years has been studied in many settings, e.g., \citep{jain2014iterative,gunasekar2015unified,elenberg2018restricted}.

We will assume the following condition, for some constants $\eps\geq 0$, 
$\alpha_f,\alpha_g\geq0$, and $c_f,c_g\in(0,+\infty]$, and
some  positive definite matrix $\Sigma\succ 0$: \begin{assumption}[Restricted Strong Convexity]\label{asm:rsc}
There exists a feasible point $(\xh,\yh)$ and subgradients  $\xi_{\yh}\in\partial f(\xh),  \zeta_{\yh}\in\partial g(\yh)$,
such that
\begin{multline}\label{eqn:rsc}\biginner{\left(\begin{array}{c}x -\xh\\ y -\yh \end{array}\right)}{\left(\begin{array}{c} \xi_x-\xi_{\xh} \\ \zeta_y -\zeta_{\yh}\end{array}\right)} \\{}\geq \alpha_f\min\{\norm{x-\xh}^2_2,c_f\norm{x-\xh}_2\} +  \alpha_g\min\{\norm{y-\yh}^2_2,c_g\norm{y-\yh}_2\}\\- \frac{1}{2}\norm{Ax + By - c}^2_\Sigma - \eps^2,\end{multline}
 for all $x\in\dom(f)$, $y\in\dom(g)$, $\xi_x\in\partial f(x)$, and $\zeta_y\in\partial g(y)$.
 \end{assumption}
 
\paragraph{Motivation} To motivate this condition, consider a first-order optimal point $(\xh,\yh)$.
We first observe that if 
the functions $f$ and $g$ were $\alpha_f$-strongly convex and $\alpha_g$-strongly convex, respectively, then we would 
have
\[\biginner{\left(\begin{array}{c}x - \xh\\ y - \yh\end{array}\right)}{\left(\begin{array}{c}\xi_x-\xi_{\xh} \\ \zeta_y -\zeta_{\yh}\end{array}\right)} 
\geq \alpha_f\norm{x-\xh}^2_2 + \alpha_g\norm{y-\yh}^2_2 \ \forall x\in\dom(f),y\in\dom(g),\, \forall \xi_x,\zeta_y.\]
If instead $f$ and/or $g$ does not satisfy strong convexity (or may even be nonconvex)
but strong convexity is regained once we impose the constraint $Ax+By=c$, we might instead have a bound of the form
\[\biginner{\left(\begin{array}{c}x - \xh\\ y - \yh\end{array}\right)}{\left(\begin{array}{c}\xi_x-\xi_{\xh} \\ \zeta_y -\zeta_{\yh}\end{array}\right)} 
\geq \alpha_f\norm{x-\xh}^2_2 + \alpha_g\norm{y-\yh}^2_2 \ \forall \text{ feasible }(x,y),\ \forall \xi_x,\zeta_y.\]
This is strictly weaker than requiring $f$ and $g$ to each be strongly convex;
here, the requirement of strong convexity is restricted to the subspace defined by the constraint $Ax+By=c$.

To accommodate the setting of ADMM, where the constraint $Ax+By=c$ is not satisfied
exactly at finite iterations, we will need to extend the statement above to allow for points that violate this constraint.
This is the motivation for subtracting the term $ \frac{1}{2}\norm{Ax + By - c}^2_\Sigma$ on the right-hand side
of~\eqref{eqn:rsc}, which
allows the strong convexity requirement to be relaxed outside of the subspace where the constraint holds.
Finally, the additional term $\eps^2$ subtracted on the right-hand side is typically a very small positive
constant, allowing for 
minor violations of the RSC property---we will return to the meaning and interpretation of this term below.

\paragraph{Parameters for the RSC condition}
We next examine the choices of constants $\alpha_f,\alpha_g,c_f,c_g$, the penalty matrix $\Sigma$,
and the ``tolerance'' term $\eps$, in this condition.
\begin{itemize}
\item {\bf Constants $\alpha_f,\alpha_g,c_f,c_g$.}
As seen earlier, in some cases the objective function may offer strong convexity in  feasible directions (i.e., $(x,y)$ such that $Ax+By=c$).
In such a case, we would take $c_f=c_g=+\infty$ (and $\eps=0$).
In other settings, however, it may not be possible to guarantee this type of strong curvature, but
we can ensure a weaker property by taking finite $c_f,c_g$. This would arise if, e.g., $f$ is a logistic loss function, which is convex globally but
is strongly convex only locally; moreover,
 in Section~\ref{sec:examples_rsc}, we will also see this type of weaker convexity guarantee
 for a sparse quantile regression problem. It may also be the case that the objective function offers strong
 convexity in the $x$ direction but may not be strongly convex in the $y$ direction (or vice versa), in which case we might have $\alpha_f>0$
 but $\alpha_g=0$, for example.
 \item {\bf Penalty matrix $\Sigma$.} The matrix $\Sigma$ appears in both the RSC assumption and in the ADMM algorithm,
 where it
enforces the constraint $Ax+By=c$. 
In other words, our assumption is that RSC holds with the same matrix $\Sigma$ as the one used in ADMM.
The RSC property therefore provides some insight into the role of the ADMM step size parameter. We can see that, in the presence of nonconvexity---or even if the problem is convex, but not globally strongly convex---the RSC property
may fail if the ADMM parameter $\Sigma$ is chosen to be too small.

While for specific problems we may have theoretical results that guide our choice of $\Sigma$ (as for
the quantile regression example---see Section~\ref{sec:examples_rsc}), more generally in practice we may
need to tune $\Sigma$ to achieve good convergence of ADMM. It is common to choose a multiple of the identity, i.e., $\Sigma=\sigma\ident_k$,
so that we only have a single scalar parameter $\sigma>0$ to tune. (In the ADMM literature,
this parameter is typically denoted by $\rho$.) In our theory, we allow
for a general $\Sigma$ rather than requiring a multiple of the identity, since in certain settings
it may be advantageous to choose a different form for $\Sigma$; we will see an 
example of this in the CT imaging application, in Section~\ref{sec:admm_for_CT}.

 \item {\bf Tolerance level $\eps$.} Finally we discuss the role of the scalar $\eps\geq 0$. This parameter allows for the condition to hold up to a small tolerance level,
and is typically taken to be vanishing, or even zero. We will see in our theoretical convergence guarantee below,
that the RSC property with a nonzero $\eps$ only guarantees convergence to within distance $\asymp\eps$ of $(\xh,\yh)$.

For example, if the optimization problem arises from a statistical question where we would like to estimate some true distribution parameters based
on a sample of size $n$, then often the function $f$ or $g$
reflects an empirical loss that is a random perturbation of some underlying ``true'' loss function.
Allowing for $\eps\asymp n^{-1/2}$ means that the RSC property
can hold even if the strong convexity properties of the underlying true loss are not preserved
exactly by the empirical loss.
The fact that the RSC property only guarantees convergence to within distance $\eps$ of the true parameters, is not worrisome in this statistical setting, because
convergence beyond the accuracy level $\eps\asymp n^{-1/2}$ is not informative---this is because a sample of size $n$ can only recover parameters up to errors of 
order $n^{-1/2}$ even with limitless computational resources (see, e.g., \citet[Section 4.1]{loh2015regularized} for further
discussion of the role of the $\eps$ term in RSC type results for high-dimensional statistics).
As an example, the scaling $\eps\asymp n^{-1/2}$ arises in the sparse quantile
regression application, for which the RSC property is studied in Section~\ref{sec:examples_rsc}.
 \end{itemize}

In Appendix~\ref{app:rsc}, we give some additional intuition and interpretations for the RSC property, for
the $\Sigma$ in particular, showing how RSC relates
to the convexity of the augmented Lagrangian $\Lcal_{\Sigma}$ defined in~\eqref{eqn:augLagr}.

\subsection{First-order conditions}\label{sec:firstorder}

A first-order stationary point (FOSP) of the optimization problem
is a feasible point $(x,y)$ such that, for any  feasible $(x',y')$, it holds that
\begin{equation}\label{eqn:firstorder}\biginner{\left(\begin{array}{c}x'-x\\y'-y\end{array}\right)}{\left(\begin{array}{c}\xi_x\\\zeta_y\end{array}\right)}\geq 0\end{equation}
for some $\xi_x\in\partial f(x)$ and some $\zeta_y\in\partial g(y)$.  In particular, for any triple $(x,y,u)\in\dom(f)\times\dom(g)\times\R^k$,
if it holds that
\begin{equation}\label{eqn:firstorder_triple}\begin{cases}
Ax + By =c,\\
- A^\top u \in \partial f(x) ,\\
- B^\top u \in \partial g(y),
\end{cases}\end{equation}
then we can verify that $(x,y)$ is a FOSP (by taking $\xi_x = -A^\top u$ and $\zeta_y = -B^\top u$ in~\eqref{eqn:firstorder}).

To prove (approximate) convergence to the target $(\xh,\yh)$, we will need to assume that this point is (approximately) first-order
optimal. 
\begin{assumption}\label{asm:approx_firstorder}
For some $\eps_{\textnormal{FOSP}}\geq 0$, 
the point $(\xh,\yh)$ satisfies 
\begin{equation}\label{eqn:approx_firstorder}\begin{cases}
A\xh + B\yh =c,\\
\norm{- A^\top \uh - \xi_{\xh}}_2 \leq   \min\left\{\frac{\alpha_fc_f}{2} , \sqrt{\alpha_f}\cdot\eps_{\textnormal{FOSP}}\right\},\\
\norm{- B^\top \uh -\zeta_{\yh}}_2 \leq  \min\left\{\frac{\alpha_gc_g}{2}, \sqrt{\alpha_g}\cdot\eps_{\textnormal{FOSP}}\right\},
\end{cases}\end{equation}
for some $\uh\in\R^k$, 
where constants $\alpha_f,\alpha_g,c_f,c_g$ 
and subgradients $\xi_{\xh},\zeta_{\yh}$
are the same as the ones appearing in Assumption~\ref{asm:rsc}.
\end{assumption}
\noindent For intuition, we can see that if $(\xh,\yh,\uh)$ were to satisfy the conditions~\eqref{eqn:firstorder_triple} exactly,
then this assumption would hold with $\eps_{\textnormal{FOSP}}=0$.

 Analogous to the role of $\eps$ in the restricted strong convexity condition, here $\eps_{\textnormal{FOSP}}$ is a tolerance
level, allowing the first-order optimality conditions to hold only approximately.
We will see that convergence is then guaranteed only up to an accuracy level that scales with these tolerance parameters
$\eps$ and $\eps_{\textnormal{FOSP}}$. 
 
A key motivation can again be found by considering a statistical setting, where
we are minimizing a loss derived from a finite sample of size $n$ (e.g., empirical risk minimization),
then we would expect the true parameters $(\xh,\yh)$ to be approximately first-order optimal
with $\eps_{\textnormal{FOSP}}\asymp n^{-1/2}$,
reflecting the usual error rates obtained with a sample size $n$.

\subsection{Main result: convergence guarantee}

Our main result proves that the ADMM iterates $(x_t,y_t,u_t)$  converge to $(\xh,\yh,\uh)$ (up to a tolerance level
determined by $\eps$ and $\eps_{\textnormal{FOSP}}$), as long as we choose 
the step size matrices $H_f,H_g$ to satisfy
\begin{equation}\label{eqn:stepsize}\begin{cases}H_f\succeq 0, \ H_f + A^\top\Sigma A\succ 0,\textnormal{ and }H_f\succeq \nabla^2 f_d(x)\textnormal{ for all $x\in\dom(f)$},\\
H_g\succeq 0,\ H_g + B^\top \Sigma B \succ 0,\textnormal{ and }H_g\succeq \nabla^2 g_d(y)\textnormal{ for all $y\in\dom(g)$}.\end{cases}\end{equation}
We note that, if $f_d$ (respectively $g_d$) is concave and $A^\top\Sigma A$ (respectively $B^\top\Sigma B$) is full-rank, then the corresponding step size matrix $H_f$ (respectively $H_g$), can be chosen to be zero.
However, even in such a setting, we may prefer to take a nonzero step size matrix for easier update step calculations, as discussed above.
We can also observe that the condition $H_f + A^\top\Sigma A\succ 0$, together with the assumption that $f_c$ is convex, proper, and lower semi-continuous, ensures that 
$x_{t+1}$ is unique and well-defined (i.e., the subproblem for the $x$ update step has a unique minimum), and similarly the condition $H_g + B^\top\Sigma B\succ 0$
ensures the same for the $y$ update step.

\begin{theorem}\label{thm1}
Suppose that the point $(\xh,\yh)$ is feasible,  satisfies
 Assumption~\ref{asm:rsc} (restricted strong convexity), and satisfies
 Assumption~\ref{asm:approx_firstorder} (approximate first-order optimality) for some $\uh\in\R^k$.
Suppose that the  nonconvex ADMM algorithm given in Algorithm~\ref{algo:main} is run with the penalty matrix $\Sigma$
chosen according to  the restricted strong convexity property~\eqref{eqn:rsc},
with step size matrices $H_f,H_g$ satisfying~\eqref{eqn:stepsize}, and  initialized at an arbitrary point $(x_0,y_0,u_0)\in\dom(f)\times\dom(g)\times\R^k$.

Define
\[\bar{x}_T = \frac{1}{T}\sum_{t=1}^{T}x_t \textnormal{ and }\bar{y}_T = \frac{1}{T}\sum_{t=1}^T y_t,\]
where $x_t,y_t$ are the iterates of the nonconvex ADMM algorithm.
Then for all $T\geq 1$, 
\begin{multline*} \alpha_f\min\left\{\norm{\bar{x}_T - \xh}^2_2,c_f\norm{\bar{x}_T-\xh}_2\right\}
+ \alpha_g\min\left\{\norm{\bar{y}_T - \yh}^2_2,c_g\norm{\bar{y}_T-\yh}_2\right\}\\
\leq  \frac{C(\xh,\yh,\uh;x_0,y_0,u_0)}{T} + 4(\eps^2 + \eps_{\textnormal{FOSP}}^2).\end{multline*}
\end{theorem}
\noindent
The function $C$ appearing in the upper bound is defined explicitly in the proof, and does not depend on the iteration number $T$.

An important observation is that convergence is guaranteed only up to the error level scaling as $\eps^2+\eps_{\textnormal{FOSP}}^2$---these
terms do not vanish
as $T\rightarrow\infty$.  To understand why this is exactly as expected, we can again consider a statistical setting, where the true parameters $(\xh,\yh)$ are estimated by minimizing a loss derived from a finite sample of size $n$; in this type of setting, convergence
can only be expected to recover $(\xh,\yh)$ up to some accuracy level. Indeed, even if we were
able to compute the global minimizer of the optimization problem, we would still expect nonzero error in recovering $(\xh,\yh)$.
In particular, as described above, in such settings we expect the RSC property and the approximate first-order optimality property
to hold with $\eps,\eps_{\textnormal{FOSP}}\asymp n^{-1/2}$; this then implies that, for sufficiently large $T$,
we have $\norm{\bar{x}_T - \xh}_2 \lesssim n^{-1/2}$. 
As discussed earlier, since this is the expected rate for parameter estimation based on a sample of size $n$
(in particular, even the {\em global} minimizer of the optimization problem will have this same error rate), we cannot
hope for a better guarantee.

\paragraph{Comparison to related work} In Section~\ref{sec:prior_work_nonconvex}, 
we discussed prior work on different variants of the nonconvex ADMM algorithm (with or without linear approximations
to the differentiable components $f_d$ and $g_d$ of the objective function). These existing results all require that at least one of the two functions ($f$ or $g$)
must be smooth, or alternatively proves a weaker convergence result, establishing properties of the limit point under the assumption that the algorithm converges 
(without proving that convergence must occur). The related MOCCA algorithm, discussed in Section~\ref{sec:prior_work_mocca}, does allow
for both $f$ and $g$ to be nonsmooth, but the convergence guarantee comes at the cost of an ``inner loop'' in the algorithm that increases in length with every iteration,
which would be extremely inefficient in practice. The contribution of Theorem~\ref{thm1} is that we can be assured that,
with the RSC assumption, the nonconvex ADMM
algorithm will converge even when both  $f$ and $g$ are nonsmooth.

\subsection{Proof of Theorem~\ref{thm1}}
Fix any point $(x,y,u)$ satisfying $Ax+By=c$. In Appendix~\ref{app:complete_proof_thm1}, we will prove that the assumption~\eqref{eqn:stepsize} on the step size matrices $H_f,H_g$
ensures that, for all $T\geq 1$, there exist some $\xi_{x_2},\dots,\xi_{x_{T+1}}$ and some $\zeta_{y_1},\dots,\zeta_{y_T}$ such that
\begin{multline}\label{eqn:mainstep}
\sum_{t=0}^{T-1} \biginner{\left(\begin{array}{c}x_{t+2}-x\\ y_{t+1}-y\end{array}\right)}{\left(\begin{array}{c}\xi_{x_{t+2}} +A^\top u\\\zeta_{y_{t+1}} +B^\top u\end{array}\right)}
+ \frac{1}{2}\sum_{t=0}^{T-1}\norm{Ax_{t+2}+By_{t+1}-c}^2_\Sigma\\
{}\leq  C_1(x,y,u;x_0,y_0,u_0) .
\end{multline}
The function $ C_1$ will be defined in the Appendix (see~\eqref{eqn:mainstep_define_C}).

Moreover, applying the restricted strong convexity assumption (Assumption~\ref{asm:rsc}), we have
\begin{multline}\label{eqn:use_rsc_in_proof}\biginner{\left(\begin{array}{c}x_{t+2}-\xh\\ y_{t+1}-\yh\end{array}\right)}{\left(\begin{array}{c}\xi_{x_{t+2}} -\xi_{\xh}\\\zeta_{y_{t+1}} -\zeta_{\yh}\end{array}\right)}
\\
{}\geq \alpha_f\min\{\norm{x_{t+2}-\xh}^2_2,c_f\norm{x_{t+2}-\xh}_2\} +\alpha_g\min\{\norm{y_{t+1}-\yh}^2_2,c_g\norm{y_{t+1}-\yh}_2\} \\- \frac{1}{2}\norm{Ax_{t+2} + By_{t+1} - c}^2_\Sigma - \eps^2\end{multline}
for each $t=0,\dots,T-1$.

Combining all of these calculations with the bound~\eqref{eqn:mainstep} above applied to $(x,y,u) = (\xh,\yh,\uh)$, and rearranging terms, we obtain
\begin{multline*}
\sum_{t=0}^{T-1} \left(  \alpha_f\min\{\norm{x_{t+2}-\xh}^2_2,c_f\norm{x_{t+2}-\xh}_2\} +\alpha_g\min\{\norm{y_{t+1}-\yh}^2_2,c_g\norm{y_{t+1}-\yh}_2\})\right)
\\\leq \sum_{t=0}^{T-1} \biginner{\left(\begin{array}{c}x_{t+2}-\xh\\ y_{t+1}-\yh\end{array}\right)}{\left(\begin{array}{c}-A^\top\uh - \xi_{\xh}\\-B^\top\uh -\zeta_{\yh}\end{array}\right)}+ C_1(\xh,\yh,\uh;x_0,y_0,u_0) + T\eps^2.\end{multline*}
Next for each $t$, we apply Assumption~\ref{asm:approx_firstorder} to calculate
\begin{multline*}\inner{x_{t+2}-\xh}{-A^\top\uh - \xi_{\xh}}\leq \norm{x_{t+2}-\xh}_2\cdot \norm{-A^\top\uh- \xi_{\xh}}_2 \\
 \leq  \min\left\{\frac{\alpha_fc_f}{2} \cdot \norm{x_{t+2}-\xh}_2, \sqrt{\alpha_f}\cdot\eps_{\textnormal{FOSP}}\cdot \norm{x_{t+2}-\xh}_2\right\}\\
 \leq  \min\left\{\frac{\alpha_fc_f}{2} \cdot \norm{x_{t+2}-\xh}_2, \frac{\alpha_f}{2}\norm{x_{t+2}-\xh}^2_2 + \frac{\eps_{\textnormal{FOSP}}^2}{2}\right\},
 \end{multline*}
and similarly for the $y$ term. Therefore,
we can rearrange the above to
\begin{multline*}
\sum_{t=0}^{T-1} \left(  \frac{\alpha_f}{2}\min\{\norm{x_{t+2}-\xh}^2_2,c_f\norm{x_{t+2}-\xh}_2\} +\frac{\alpha_g}{2}\min\{\norm{y_{t+1}-\yh}^2_2,c_g\norm{y_{t+1}-\yh}_2\})\right)
\\\leq C_1(\xh,\yh,\uh;x_0,y_0,u_0) + T\eps^2 + T\eps_{\textnormal{FOSP}}^2.\end{multline*}
Next,  noting that $x_1$ is a deterministic function of $(x_0,y_0,u_0)$, we  define \[C(x,y,u;x_0,y_0,u_0) = 4C_1(x,y,u;x_0,y_0,u_0) +  2\alpha_f\min\{\norm{x_1-\xh}^2_2,c_f\norm{x_1-\xh}_2\} .\]
We can then relax the bound above to
\begin{multline}\label{eqn:thm1_laststep}
\sum_{t=0}^{T-1} \left(  \frac{\alpha_f}{2}\min\{\norm{x_{t+1}-\xh}^2_2,c_f\norm{x_{t+1}-\xh}_2\} +\frac{\alpha_g}{2}\min\{\norm{y_{t+1}-\yh}^2_2,c_g\norm{y_{t+1}-\yh}_2\})\right)
\\\leq \frac{1}{4}C(\xh,\yh,\uh;x_0,y_0,u_0) + T\eps^2 + T\eps_{\textnormal{FOSP}}^2.\end{multline}
Next we will use the following elementary fact: for any nonnegative $c,r_1,\dots,r_n$,
\[\sum_{i=1}^n \min\{r_i^2,cr_i\} \geq \frac{n}{2}\min\left\{\left(\frac{1}{n}\sum_{i=1}^n r_i\right)^2,c\left(\frac{1}{n}\sum_{i=1}^n r_i\right)\right\}.\]
Therefore, applying this with $\norm{x_{t+1}-\xh}_2$ in place of the $r_i$ terms, we have
\[\sum_{t=0}^{T-1} \frac{\alpha_f}{2}\min\{\norm{x_{t+1}-\xh}^2_2,c_f\norm{x_{t+1}-\xh}_2\} 
\geq \frac{T\alpha_f}{4}\min\left\{\norm{\bar{x}_T - \xh}^2_2,c_f\norm{\bar{x}_T-\xh}_2\right\},\]
where the last step holds since $\frac{1}{T}\sum_{t=0}^{T-1}\norm{x_{t+1}-\xh}_2\geq \norm{\bar{x}_T-\xh}_2$ by convexity.
An analogous bound holds for the $y$ term. Combining this with~\eqref{eqn:thm1_laststep} completes the proof.

\section{Example: sparse high-dimensional quantile regression}\label{sec:examples}
In this section, we will develop a concrete example of our framework,
to illustrate the empirical performance and convergence properties of our method. Consider a regression setting where
\[w_i = \phi_i^\top \xh + \textnormal{(noise)},\quad i=1,\dots,n,\]
for a sparse true signal $\xh\in\R^d$. The response variables $w_i\in\R$ and the sensing matrix $\Phi = (\phi_1,\dots,\phi_n)^\top\in\R^{n\times d}$
are observed, and our goal is to recover $\xh$. If the noise is heavy-tailed, then a standard least-squares regression may perform poorly,
and we may prefer the more robust properties of a quantile regression.
Specifically, for any desired quantile $q\in(0,1)$, consider the quantile loss
\[\ell_q(t) = q \cdot\max\{t,0\} + (1-q)\cdot \max\{-t,0\}.\]
Then if we seek to minimize 
\[\frac{1}{n}\sum_{i=1}^n \ell_q\big(w_i - \phi_i^\top x\big)\]
over $x\in\R^d$, this loss corresponds to aiming for $\phi_i^\top x$ to equal the $q$-th quantile of $w_i$.
(Note that for the special case $q =0.5$, i.e., median regression, this loss is equal to the $\ell_1$ norm, up to rescaling.)

\begin{figure}[t]\centering
\includegraphics[height=0.3\textwidth]{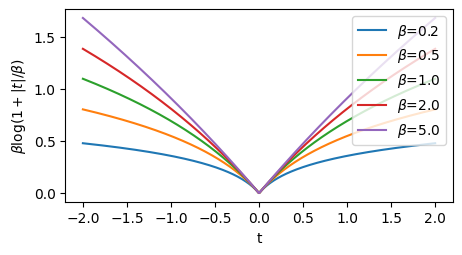}
\label{fig:logL1}\caption{Illustration of the nonconvex sparsity-promoting penalty $\sum_j \beta\log(1+|x_j|/\beta)$ that appears
in the objective function~\eqref{eqn:sparseQR} for the sparse high-dimensional quantile regression example.
 The figure plots the function $t\mapsto \beta\log(1+|t|/\beta)$, for a range of values of $\beta$. 
The functions are all nondifferentiable at $t=0$, and are similar to the absolute value function  for $t\approx 0$, but smaller values of $\beta$ correspond to greater
nonconvexity as $|t|$ increases.}
\end{figure}

In the high-dimensional setting where $n<d$, minimizing this loss is not meaningful (in general, we can always find a vector $x\in\R^d$ 
that interpolates the data, i.e., $\phi_i^\top x=w_i$ for all $i$, which clearly leads to overfitting). We will therefore consider a penalized
version of this loss:
\begin{equation}\label{eqn:sparseQR}\argmin_{x\in\R^d} \big\{\loss(x)\big\}\textnormal{ where }\loss(x)=\frac{1}{n}\sum_{i=1}^n \ell_q\big(w_i - \phi_i^\top x\big) + \lambda\sum_{j=1}^d \beta\log(1 + |x_j|/\beta).\end{equation}
The last term is a nonconvex regularizer that encourages a sparse solution; see \citet{fazel2003log,candes2008enhancing} for background. For $\beta=+\infty$, the regularizer is equal to the $\ell_1$ norm, a standard
convex penalty for recovering sparse signals, while $\beta<+\infty$ leads to a nonconvex penalty. 
Smaller values of $\beta$ correspond to greater nonconvexity, which makes the optimization problem more challenging
but comes with the benefit of less shrinkage on the nonzero values in the signal vector $x$ (see Figure~\ref{fig:logL1}).

To enable theoretical guarantees, we will add one small modification to this optimization problem, and will instead solve
\begin{equation}\label{eqn:sparseQR_ball}\argmin_{x\in\R^d : \norm{x}_2\leq R} \big\{\loss(x)\big\}\end{equation}
for a large radius $R$, where this constraint
is added to ensure that the iterations $x$ do not diverge to infinity. We will see in our theoretical results that we 
can set $R$ to be extremely large without compromising the convergence guarantee; in practice, therefore, we would expect that
iteratively solving~\eqref{eqn:sparseQR_ball} would be indistinguishable from iteratively solving the unconstrained version~\eqref{eqn:sparseQR},
since the constraint $\norm{x}_2\leq R$ would likely never be active.

\subsection{Implementing nonconvex ADMM}\label{sec:examples_implement_admm}
For the sparse quantile regression problem~\eqref{eqn:sparseQR_ball}, we will introduce an additional variable $y$ (with the constraint $y=\Phi x$) so that the 
optimization problem can be solved with Algorithm~\ref{algo:main}---we will minimize
\[\argmin_{x\in\R^d,y\in\R^n} \left\{\frac{1}{n}\sum_{i=1}^n \ell_q\big(w_i - y_i \big) + \lambda\sum_{j=1}^d \beta\log(1 + |x_j|/\beta) \ : \ y = \Phi x,\norm{x}_2\leq R\right\}.\]
To solve~\eqref{eqn:sparseQR_ball}, we define $A=\Phi$, $B = - \ident_n$, and $c=0$, and 
 run  Algorithm~\ref{algo:main} with  parameters $\Sigma=\sigma\ident_n$ (for a chosen value of the tuning parameter $\sigma>0$), $H_f=\sigma(\gamma\ident_d - \Phi^\top \Phi)$ (with $\gamma = \norm{\Phi}^2$ so that $H_f\succeq 0$), and $H_g=0$,
and with functions 
\[f_c(x) = \lambda\norm{x}_1 + \delta_{\norm{x}_2\leq R},\quad f_d(x) =   \lambda\sum_{j=1}^d \left(\beta\log(1 + |x_j|/\beta) - |x_j|\right),\]
where $ \delta_{\norm{x}_2\leq R}$ is the convex indicator function (i.e., $ \delta_{\norm{x}_2\leq R}=0$ if $\norm{x}_2\leq R$, and $ \delta_{\norm{x}_2\leq R}=+\infty$
otherwise),
and with
\[ g_c(y) = \frac{1}{n}\sum_{i=1}^n \ell_q\big(w_i - y_i \big) ,\quad  g_d(y) \equiv 0.\]
The update steps for Algorithm~\ref{algo:main} can  be calculated in closed form (details are given in Appendix~\ref{app:sparseQR}).
We note that the function $f_d$ is concave and twice differentiable, with $\nabla^2 f_d(x)\succeq -\lambda\beta^{-1}\ident_d$ for all $x$, so its concavity is bounded.

\subsection{Theoretical results}\label{sec:examples_rsc}
Our theoretical results guarantee convergence for the nonconvex ADMM algorithm
as long as the RSC property~\eqref{eqn:rsc} and the the approximate first-order optimality property~\eqref{eqn:approx_firstorder}
both hold, to verify the assumptions of Theorem~\ref{thm1}. In particular, RSC-type properties for sparse
high-dimensional quantile regression have been studied in the literature, e.g., see \citet[Lemma C.3]{zhao2014general} or \citet[Lemma 4]{belloni2011L1}. The conditions proved in the literature 
appear in a different form than the RSC property studied here, so we verify that the property~\eqref{eqn:rsc} holds under some mild assumptions.
The following result is proved in Appendix~\ref{app:sparseQR_proof}.
\begin{proposition}\label{prop:quantile_rsc}
Suppose that the observations are given by
\[w_i = \phi_i^\top \xh + z_i, \quad i=1,\dots,n\]
for some sample size $n\geq4$, and let $\yh=\Phi\xh$. Assume that:
\begin{itemize}
\item The feature vectors $\phi_i\in\R^d$ are i.i.d.~with distribution $\mathcal{D}_\phi$,
where for $\phi\sim \mathcal{D}_\phi$, it holds that $\norm{\phi}_{\infty}\leq B_\phi$ almost surely,
and that $\EE{|\phi^\top u|^2}\geq a_\phi$ and $\EE{|\phi^\top u|^3}\leq b_\phi$ for any fixed unit vector $u\in\R^d$;
\item The noise terms $z_i\in\R$ are drawn independently from the feature vectors $\phi_i$,
and moreover are  i.i.d.~with density $h_z$, for which $z=0$ is the $q$-th quantile,
and which satisfies  $h_z(t)\geq c_z$ for all $|t|\leq t_z$, for some $c_z,t_z>0$;
\item The true vector $\xh$ has at most $s_*$ nonzero entries, where
\[1\leq s_* \leq C_0 \cdot \frac{n}{\log(nd)} \]
for a constant $C_0>0$ that depends only on $c_z,t_z,a_\phi,b_\phi,B_\phi$;
\item  The parameters $\lambda,\beta,R$ are chosen to satisfy
\[\lambda = C_\lambda \sqrt{\frac{\log(nd)}{n}}\textnormal{ for some }C_\lambda \in \left[C_1, C_1 \sqrt{\frac{C_0 \cdot \frac{n}{\log(nd)}}{s_*}}\, \right]\]
and
\[ R\geq \norm{\xh}_2\text{  \ and \  }\beta \geq C_\lambda \max\{1,R\}\cdot C_2 \sqrt{\frac{\log(nd)}{n}},\]
for constants $C_1,C_2>0$ that depend only on $c_z,t_z,a_\phi,b_\phi,B_\phi$.
\end{itemize}

Then, for any $\sigma>0$, with probability at least $1-(nd)^{-1}$, the RSC property~\eqref{eqn:rsc} holds with  
\[\alpha_f = C_3,\quad \alpha_g= 0, \quad c_f = c_g=1,\quad\Sigma=\sigma\ident_n, \quad \eps^2 = C_4\max\{1,\sigma^{-1}\}\cdot \frac{s_*\log(nd)}{n},\]
and the approximate first-order optimality property~\eqref{eqn:approx_firstorder} holds with
\[\eps_{\textnormal{FOSP}}^2 = C_5\cdot \frac{s_*\log(nd)}{n},\]
where $C_3,C_4,C_5>0$ are constants that depend only on $c_z,t_z,a_\phi,b_\phi,B_\phi$ and on $C_\lambda$.
\end{proposition}
With this result in place, if $\lambda,\beta,R$ are chosen appropriately,
then Theorem~\ref{thm1} ensures that, after $T$ iterations of ADMM, the estimate $\bar{x}_T$ will satisfy
\[\min\{\norm{\bar{x}_T - \xh}^2_2,\norm{\bar{x}_T-\xh}_2\} \leq \O{ \frac{1}{T} + \frac{s_*\log(nd)}{n}},\]
which we can simplify to
\[\norm{\bar{x}_T - \xh}_2\leq \O{\sqrt{\frac{1}{T}+\frac{s_*\log(nd)}{n}}}.\]
In contrast, the minimax error rate for estimating $\xh$, in this high-dimensional sparse
regression setting, is $\O{\sqrt{\frac{s_*\log(d/s_*)}{n}}}$ \citep[Theorem 1(b)]{raskutti2011minimax}.
This shows that, up to a slightly different log factor, the error of $\bar{x}_T$
matches the minimax rate once $T$ is sufficiently large.

\paragraph{Comparing to existing theory}
As discussed in Section~\ref{sec:prior_work_nonconvex},
previous results  establishing convergence for nonconvex ADMM assume, at minimum, that either $f$ or $g$
is differentiable and has a Lipschitz gradient. We can see immediately that this property is violated for the 
sparse quantile regression problem~\eqref{eqn:sparseQR} (or for its constrained version~\eqref{eqn:sparseQR_ball}), since the functions $f$ and $g$ are both nondifferentiable.
In contrast, our new RSC-based framework is  able to provide a guarantee,
and so
this example illustrates the flexibility and broad applicability of RSC type assumptions, as compared to other assumptions in the literature.

\begin{figure}[p]\centering
\includegraphics[width=0.49\textwidth]{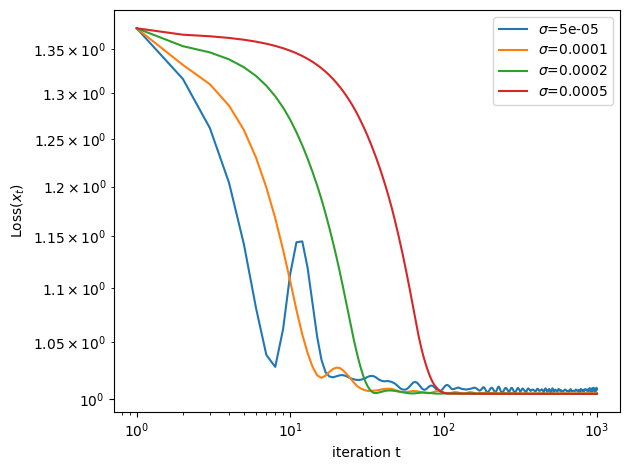}
\includegraphics[width=0.49\textwidth]{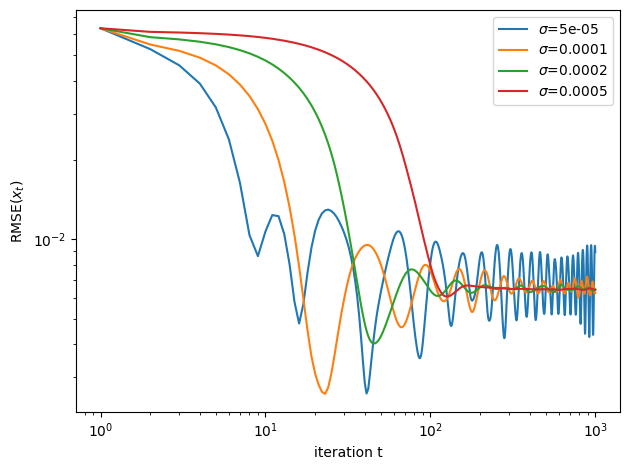}\\\medskip
 \includegraphics[width=0.49\textwidth]{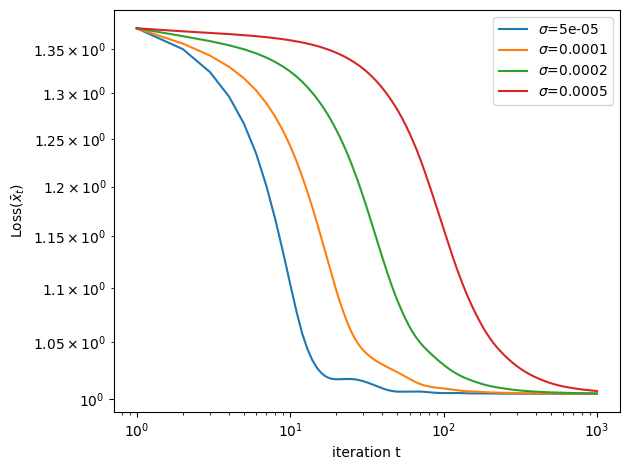}
 \includegraphics[width=0.49\textwidth]{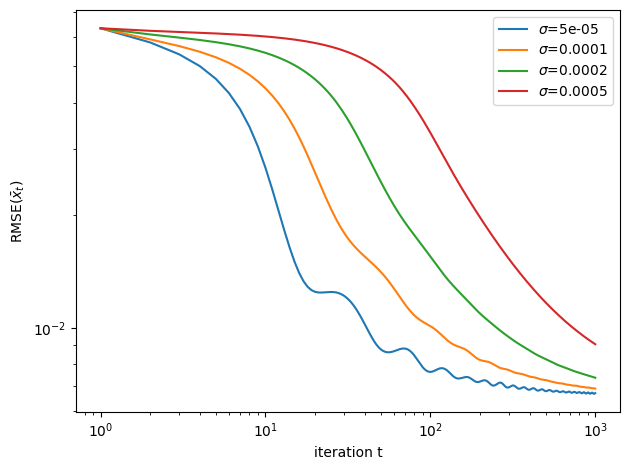}\\
\caption{Results for the sparse quantile regression example (see Section~\ref{sec:examples_empirical}). 
The figure shows the value of the objective function~\eqref{eqn:sparseQR} over iteration $t=1,\dots,500$ of the algorithm,
run with various values of the parameter $\sigma$ as shown. The top row shows the loss function value for $x_t$ (the estimate at time $t$),
as well as its root-mean-square-error (RMSE) $\frac{1}{\sqrt{d}}\|x_t - \xh\|_2$, while the bottom plot
shows the loss and the RMSE for $\bar{x}_t$ (the running average). All axes are on the log scale.}
\label{fig:examples}
\end{figure}

\subsection{Empirical results}\label{sec:examples_empirical}
We next demonstrate the performance of our algorithm on the sparse quantile regression problem.
Code reproducing the simulation and all figures is available at \url{https://github.com/rinafb/ADMM_CT}.

We choose dimension $d=2500$ and sample size $n=2000$ for a challenging high-dimensional setting. The matrix $\Phi\in\R^{n\times d}$ is constructed
with i.i.d.~$\normal(0,1)$ entries. We define
\[w_i = \phi_i^\top \xh + z_i,\]
where $\phi_i$ is the $i$th row of $\Phi$, and the true signal is given by $\xh = (1,\dots,1,0,\dots,0)$, with $s_*=10$ nonzero entries.
 The noise
terms $z_i$ are drawn i.i.d.~from $t_5$, the standard t distribution with 5 degrees of freedom, which is a heavy-tailed distribution.
We choose the quantile $q=0.5$ (i.e., a median regression).
For the penalty term, we choose $\lambda = 0.1$ and $\beta=0.5$;  this small value of $\beta$ 
means that the penalty has substantial nonconvexity (see Figure~\ref{fig:logL1}).
The parameter $\sigma$ controlling the enforcement of the constraint in ADMM (i.e., with $\Sigma = \sigma\mathbf{I}_d$ in Algorithm~\ref{algo:main})
is varied as $\sigma\in\{0.00005,0.0001,0.0002,0.0005\}$.

The results after running Algorithm~\ref{algo:main} for 1000 iterations are displayed in Figure~\ref{fig:examples}.
The plot displays the loss, $\loss(x_t)$ at each iteration $t$, where $\loss(\cdot)$ is the objective
function defined in~\eqref{eqn:sparseQR}, as well as the root-mean-square-error (RMSE), $\frac{1}{\sqrt{d}}\|x_t - \xh\|_2$. (We do not impose a constraint $\norm{x}_2\leq R$, since as mentioned above, the theory allows
for $R$ to be extremely large, and the iterations $x_t$ do not violate this constraint in practice.)
The plot also shows $\loss(\bar{x}_t)$ and $\frac{1}{\sqrt{d}}\|\bar{x}_t - \xh\|_2$, the loss and RMSE
of the running average of the estimates, $\bar{x}_t = \frac{1}{t}\sum_{t'=1}^t x_{t'}$.
 The convergence of the loss and RMSE for $\bar{x}_t$ across all $\sigma$ values
is supported by our theoretical result,  Proposition~\ref{prop:quantile_rsc}, which
shows that the RSC property holds (with high probability) for {\em any} $\sigma>0$, as long as the tolerance term $\eps$ is adjusted accordingly.
Note that the RMSE (for both $x_t$ and $\bar{x}_t$) does not converge to zero, but instead appears to be
converging to a small but positive value; this is due to the noise in the data.

Interestingly, we see that overly small values of $\sigma$ lead to some instability in the convergence of the loss and the RMSE,
suggesting that the RSC property may not be sufficient to ensure convergence of the iterates themselves (the $x_t$'s)
rather than the running averages (the $\bar{x}_t$'s).\footnote{An alternative explanation for this empirical result
is simply that the parameter $\eps$ in the RSC property~\eqref{eqn:rsc}
is larger, when $\sigma$ is chosen to be smaller, as in Proposition~\ref{prop:quantile_rsc}; since convergence is only guaranteed up to the tolerance
level $\eps$ in Theorem~\ref{thm1}, this may explain the apparent lack of convergence for $x_t$ when $\sigma$ is chosen to be very small.}
On the other hand, overly large values of $\sigma$ may lead to somewhat slower convergence; intuitively, enforcing the constraint $y=\Phi x$
with too strong of a penalty will make it difficult for the algorithm to make fast progress with alternating updates of $x$ and $y$.

%%%%%%%%%%%%%%%%%%%%%%%%%%%%%%%%%%%%%%%%%%%%%%%%%%%%%%%%%%%%%%%%%%%%%%%%%%%%%

\section{Application: CT imaging}\label{sec:CT}

We next apply our algorithm and convergence results to the problem of image reconstruction
in computed tomography (CT) imaging, which is the motivating application for this work.
In CT, we would like to reconstruct an image of an unknown object $x$ (e.g., produce a 3D image of a patient's head or abdomen, in
the setting of medical CT). The available measurements obtained from the CT scanner
consist of measuring the intensity of an X-ray beam passing through the unknown object. A lower intensity of the beam
when it reaches the detector indicates higher density in the unknown object along that ray.

We now introduce some notation to make this problem more precise. We will begin with a simple version of the problem,
and then will add additional components step by step to build intuition.
Let $x=(x_k)\in\R^{n_k}$ denote the unknown image, where
 $k=1,\dots,n_k$ indexes pixels (or voxels), after we have discretized to a 2D (or 3D) grid---for example,
in two dimensions, $n_k = N_x\cdot N_y$ for an $N_x\times N_y$ grid. 

To obtain an image, the scanner sends an X-ray beam along $n_\ell$ many rays.
For example, for many clinical scanners in a medical setting, the device rotates around  the patient, taking images from $N_{\textnormal{img}}$ 
many angles; for each of these images, there are $N_{\textnormal{cell}}$ many detector cells measuring the intensity of the beam after it passes
through the patient's body. This leads to $n_\ell = N_{\textnormal{img}}\cdot N_{\textnormal{cell}}$ many rays $\ell=1,\dots,n_\ell$ 
along which measurements are taken.

\begin{figure}\centering\begin{tabular}{@{}c@{}}
\includegraphics[width = 0.37\textwidth]{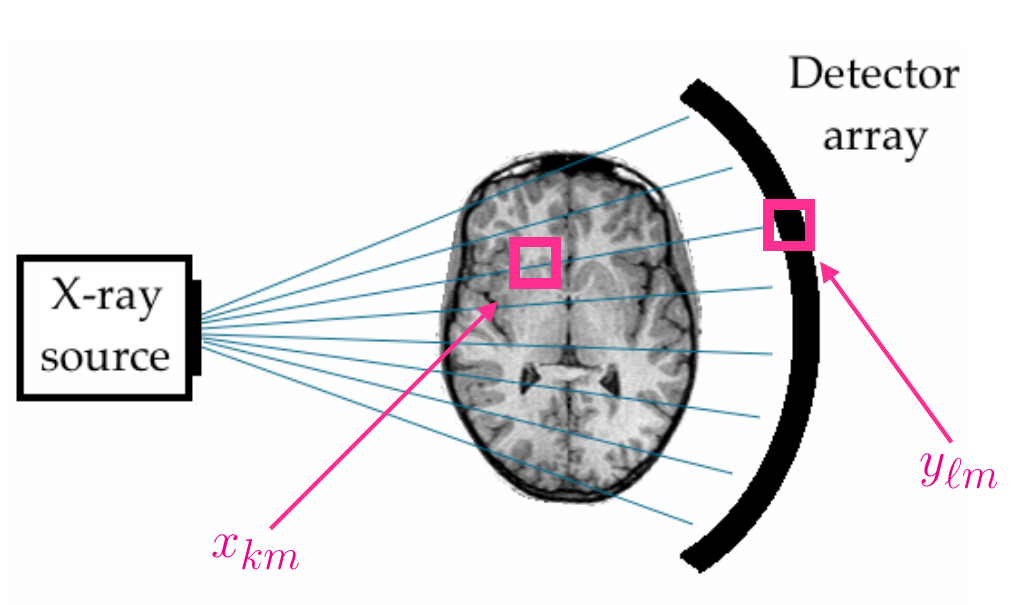}\end{tabular}
\ \begin{tabular}{@{}c@{}}
\includegraphics[width = 0.6\textwidth]{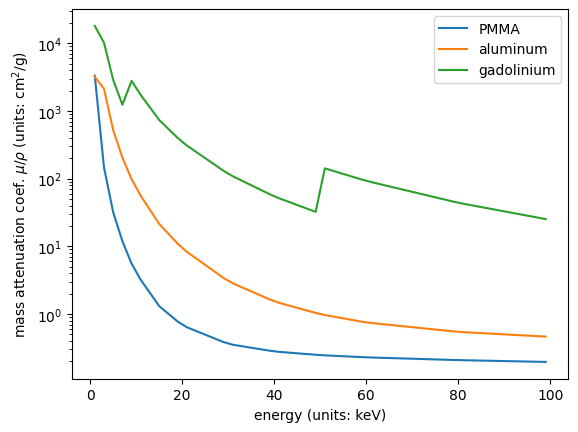}\end{tabular}
\caption{Left: schematic of the projection operator. Here $x_{k m}$ is the amount of material $m$ present at pixel $k$,
while $y_{\ell m} = (Px)_{\ell m}$ is the total amount of material $m$ present along ray $\ell$ of the scan. Right: attenuation curves for several common materials.}
\label{fig:CT_intro}
\end{figure}

Now let $P = (P_{\ell k})\in\R^{n_\ell\times n_k}$ be the projection matrix, with $P_{\ell k}$ measuring the length
of the intersection between ray $\ell$ and pixel $k$. The product $Px \in\R^{n_\ell}$ measures
the projection of the object $x$, where $(Px)_{\ell}$ measures the total amount of material that lies along ray $\ell$ (see Figure~\ref{fig:CT_intro} for a schematic). 
The attenuation (i.e., the loss of intensity) of the X-ray beam that travels along ray $\ell$ 
depends on $(Px)_\ell$. In particular, ignoring photon scattering
and other sources of noise, the measurements follow a model of the form
\[\frac{\textnormal{Intensity of the beam after passing through the object along ray $\ell$}}{\textnormal{Intensity of the beam entering the object along ray $\ell$}} \approx e^{-\mu \cdot (Px)_\ell},\]
where $\mu>0$ is called the linear attenuation coefficient. 
While most clinical scanners measure the total energy of the beam when it reaches the detector, 
here we consider a different type of hardware, {\em photon counting CT}, where the measurement is a count of the number of photons
reaching the detector. In this case, we can model this count as
\[C_\ell\sim \textnormal{Poisson}\left(S \cdot \exp\{-\mu \cdot(Px)_\ell\}\right),\]
where $S$ is the number of photons incident on the detector pixel (characterizing the intensity of the X-ray beam for a fixed time-duration scan), and $C_\ell$ 
is the number of photons reaching detector  after passing through the object along ray $\ell$. 

In fact, since different detector cells may have slightly different sensitivities, a more accurate model is
\begin{equation}\label{eqn:TPL_onematerial}C_\ell \sim \textnormal{Poisson}\left(S_\ell \cdot \exp\{-\mu \cdot(Px)_\ell\}\right),\end{equation}
where the scalar term $S_\ell$ combines beam intensity with detector sensitivity for ray $\ell$.

\paragraph{Multiple materials}
In practice, the unknown object can consist of multiple materials, which each behave differently in terms of the attenuation of the beam.
Let $m=1,\dots,n_m$ index the materials that make up the object---for example,
in a simple medical setting we might have $n_m=3$ with bone, soft tissue, and an injected contrast material such as a gadolinium
or iodine compound.
The goal is now to reconstruct the image $x = (x_{km})\in\R^{n_k\times n_m}$, 
where, for each pixel $k$, $x_{km}$ is the proportion of that pixel that is occupied by each material.
We can update our model~\eqref{eqn:TPL_onematerial} above to
\begin{equation}\label{eqn:TPL_monochromatic}C_\ell \sim \textnormal{Poisson}\left(S_\ell \cdot \exp\left\{-\sum_m\mu_m \cdot(Px)_{\ell m}\right\}\right),\end{equation}
where now $\mu_m>0$ is the (known) linear attenuation coefficient for material $m$.

\paragraph{A non-monochromatic beam}
Thus far, the Poisson model for CT image reconstruction does not introduce nonconvexity---maximizing the log-likelihood of the Poisson
model given in~\eqref{eqn:TPL_monochromatic} is a convex problem.
However, this model ignores the nature of the X-ray beam used in practice, for which the photons are distributed across a spectrum
of energies. The attenuation coefficient for a material $m$ in fact depends on the energy of the photon, with each material 
exhibiting its own attenuation curve across the range of energies---see Figure~\ref{fig:CT_intro} for an example. In particular,
in medical applications, contrast materials such as gadolinium or iodine are used for their unique attenuation curves, which make these materials
easier to distinguish from surrounding soft tissue in a CT scan.

Our model can now be updated to the following:
\begin{equation}\label{eqn:TPL_onewindow}C_\ell \sim \textnormal{Poisson}\left(\sum_i S_{\ell i} \cdot \exp\left\{-\sum_m\mu_{mi} (Px)_{\ell m}\right\}\right),\end{equation}
where $i=1,\dots,n_i$ is the index over a discretized grid of the range of energies in the X-ray beam, while $S_{\ell i}$ is the intensity
of the X-ray beam (combined with detector sensitivity) for energy level $i$ and ray $\ell$, and $\mu_{mi}$ is the attenuation coefficient for material $m$
at energy level $i$.
The photons measured by the detector may come from any energy level in the spectrum (i.e., the measurements $C_\ell$ are a combination of photons
from each energy level $i$).
The resulting log-likelihood maximization problem is no longer a convex function, which is a core challenge of CT image reconstruction.

\paragraph{Spectral CT} In spectral CT, the hardware of the scanner allows partial identification of the photon energies, making the reconstruction
problem somewhat easier. Specifically, the detectors are programmed with several thresholds, separating the range of energies 
of the beam into ``windows'' $w=1,\dots,n_w$ (for example, 2 windows in some current clinical scanners, 
or 3--5 windows in current research prototypes). The measurements are now indexed by $C_{w \ell}$, the number of photons
in energy window $w$ measured along ray $\ell$. In theory, the windows form a partition of the energy range, but in practice there
is some noise at the boundaries between windows (that is, a photon with energy near the chosen threshold has some chance of being detected in either window).
To quantify this, let $S_{w\ell i}$ 
incident photon spectral density at energy $i$, multiplied by the probability of a photon at energy $i$ being detected in window $w$ (for the detector
sensitivity corresponding to ray $\ell$). These values are typically estimated ahead of time with a calibration process. Then
the model for our measurements $C_{w \ell}$ is given by
\begin{equation}\label{eqn:TPL_spectralCT}C_{w\ell} \sim \textnormal{Poisson}\left(\sum_i S_{w\ell i} \cdot \exp\left\{-\sum_m\mu_{mi} (Px)_{\ell m}\right\}\right).\end{equation}
We can estimate the image $x$ by maximum likelihood estimation, but as before in~\eqref{eqn:TPL_onewindow}, maximizing the log-likelihood is a non-convex problem. (See \citet{barber2016algorithm} for more details on this model.)

\subsection{Image reconstruction with nonconvex ADMM}\label{sec:admm_for_CT}
We now consider the image reconstruction problem:
given observations (photon counts) $C_{w\ell}$, we would like to solve
\begin{equation}\label{eqn:CT_opt}\xh = \argmin_{x\in\R^{n_k\times n_m}}\loss(Px) ,\end{equation}
where $\loss(y)$ is the negative log-likelihood of the Poisson model for spectral CT~\eqref{eqn:TPL_spectralCT} 
given the projected object $y=Px\in\R^{n_\ell\times n_m}$:
\[\loss(y) = \sum_{w\ell} \left[\sum_i S_{w\ell i} \exp\left\{ - \sum_m\mu_{mi}y_{\ell m}\right\} - C_{w\ell}\log\left(\sum_i S_{w\ell i} \exp\left\{ - \sum_m\mu_{mi}y_{\ell m}\right\}\right)\right].\]
We note that the first term of this loss is convex in $y$ (and therefore, in $x$), while the second term is concave.

\paragraph{Modifying the exp function}
Under a well-specified model, the true image $x$ and its projection $y=Px$ must both consist of nonnegative values.
However, model misspecification, or inaccurate estimates of $x$ and/or $y$ at early stages of the iterative algorithm, can lead to negative
values. Examining the loss function, we can see that this issue may pose problems for optimization, since $t\mapsto \exp\{t\}$ has high curvature at large values
of $t$. To resolve this, we replace the $\exp\{\cdot\}$ function with the approximation:
\[\qexp\{t\} = \begin{cases}
\exp\{t\}, & t\leq 0,\\
1 + t + \frac{1}{2}t^2, & t\geq 0.\end{cases}\]
The ``q'' in the name of this modified function refers to the fact that,
 for positive values of $t$ we replace $\exp\{t\}$ with a quadratic approximation, by taking the Taylor expansion at $t=0$.
For negative values of $t$, the function is unchanged.
This choice means that the function $\qexp\{t\}$ is continuously twice differentiable and
is equal to $\exp\{t\}$ at all negative values of $t$ (i.e., for any feasible nonnegative image $x$), while at the same time ensuring 
a bounded second derivative to avoid problems in the optimization. We will therefore work with a modified loss function,
\[\loss(y) = \sum_{w\ell} \left[\sum_i S_{w\ell i} \qexp\left\{ - \sum_m\mu_{mi}y_{\ell m}\right\} - C_{w\ell}\log\left(\sum_i S_{w\ell i} \qexp\left\{ - \sum_m\mu_{mi}y_{\ell m}\right\}\right)\right].\]

It is important to note that, for CT imaging, if the model is well specified then the argument
to $\exp\{\cdot\}$ or to $\qexp\{\cdot\}$ should always be nonpositive at the true $\yh=P\xh$ 
(i.e., $\sum_m\mu_{mi}y_{\ell m}$ should be nonnegative at the true $\yh$),
and therefore, $\qexp\{\cdot\}$ should be identical to $\exp\{\cdot\}$ in the relevant range of values. Empirically,
however, the convergence behavior
of the optimization problem is often helped by allowing both positive and negative values, particularly in early iterations,
and this can also provide useful flexibility in the case of model misspecification.

\paragraph{Running nonconvex ADMM}

To reformulate the minimization problem~\eqref{eqn:CT_opt} into the setting of nonconvex ADMM,
we will solve the equivalent problem
\begin{equation}\label{eqn:CT_opt_xy}\xh,\yh = \argmin_{\substack{x\in\R^{n_k\times n_m}\\y\in\R^{n_\ell\times n_m}}}\left\{\loss(y) \ : \  Px = y\right\}.\end{equation}
Now define
$f(x)=f_c(x)=f_d(x)\equiv 0$,
and write
$g(y) = g_c(y) + g_d(y)$
where
\begin{equation}\label{eqn:CT_define_g_c}g_c(y) =  \sum_{w\ell} \sum_i S_{w\ell i} \qexp\left\{ - \sum_m\mu_{mi}y_{\ell m}\right\}\end{equation}
and
\[g_d(y) = -\sum_{w\ell}C_{w\ell} \log\left(\sum_i S_{w\ell i} \qexp\left\{ - \sum_m\mu_{mi}y_{\ell m}\right\}\right).\]
Then $\loss(y) = g(y)$, and
we have therefore reformulated the spectral CT maximum likelihood estimation problem into 
the form of our nonconvex ADMM algorithm, i.e., $\min_{x,y}\{f(x)+g(y):Px=y\}$, minimizing a composite
objective function under a linear constraint. In particular,
converting the matrix variables $x\in\R^{n_k\times n_m}$ and $y\in\R^{n_\ell\times n_m}$
to vectorized variables $\vect(x)\in\R^{n_kn_m}$ and $\vect(y)\in\R^{n_\ell n_m}$,
the constraint $Px = y$ can be rewritten as $A\vect(x)+B\vect(y)=c$ where $A = P\otimes\ident_{n_m}\in\R^{n_\ell n_m\times n_k n_m}$, $B=-\ident_{n_\ell n_m}$, and $c=0$
(here $\otimes$ denotes the matrix Kronecker product).

 We can therefore implement Algorithm~\ref{algo:main} 
for solving this optimization problem. 
To run Algorithm~\ref{algo:main} for the CT image reconstruction problem~\eqref{eqn:CT_opt_xy}, we need to choose the step size matrices $H_f,H_g$ and the penalty matrix $\Sigma$.
Following the construction proposed by \citet{pock2011diagonal} (for the convex setting),
we begin by selecting a parameter $\sigma>0$. 
We will choose step size matrix $H_g=0$ for $y$, while for the variable $x$ our step size matrix $H_f$
will be equal to
$H_f = (Q_f\otimes \ident_{n_m})- (P\otimes\ident_{n_m})^\top \cdot (\tilde\Sigma\otimes \ident_{n_m})\cdot  (P\otimes\ident_{n_m})$,
and the penalty parameter matrix $\Sigma$ will be defined as $\Sigma = \tilde\Sigma\otimes \ident_{n_m}$,
where $Q_f\in\R^{n_k\times n_k}$ and $\tilde\Sigma\in\R^{n_\ell \times n_\ell}$ are diagonal matrices with entries
\[(Q_f)_{kk} = \sigma\sum_\ell P_{\ell k}, \quad \tilde\Sigma_{\ell\ell} = \frac{\sigma}{\sum_k P_{\ell k}}.\]
With these constructions, $H_f$ is positive semidefinite as required \citep[Lemma 2]{pock2011diagonal}. 
The update steps for the nonconvex ADMM algorithm are computed in  Appendix~\ref{app:CT}.

\subsection{CT simulation}
To demonstrate the algorithm's performance on the nononvex CT image reconstruction problem,
we carry out a small-scale simulation in Python. (Performance of these methods on a large scale requires 
more careful implementation, and is addressed in our application specific work in \citet{barber2016algorithm,schmidt2020spectral}.)
Code reproducing the simulation and all figures is available at \url{https://github.com/rinafb/ADMM_CT}.

\begin{figure}\centering
\includegraphics[width = \textwidth]{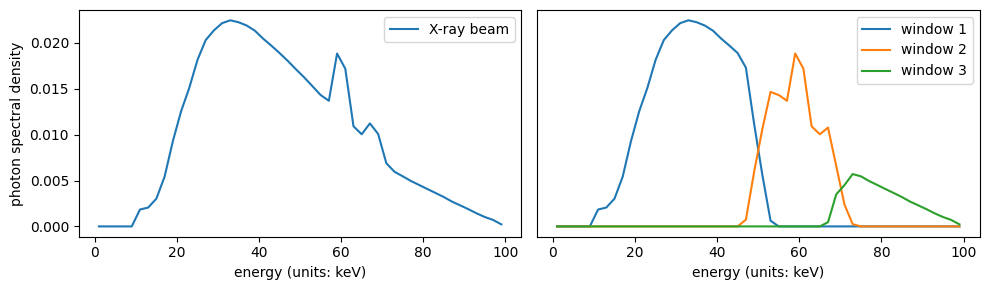}
\caption{Left: the X-ray beam spectrum. This figure displays the density of the distribution of energies in the beam, i.e., how
the total intensity of the beam is split across the energy spectrum. Right: for each energy window $w$, the displayed curve is proportional to
the spectral response parameters $S_{w\ell i}$. These values are set to be constant across all rays $\ell$,
and so the figure plots the value across all energy levels $i$ for each detector window $w$, rescaled so that the sum of the three response curves
is equal to the density plot of the X-ray beam spectrum on the left.}
\label{fig:CT_intensity}
\end{figure}

\begin{figure}\centering
\includegraphics[width = 0.7\textwidth]{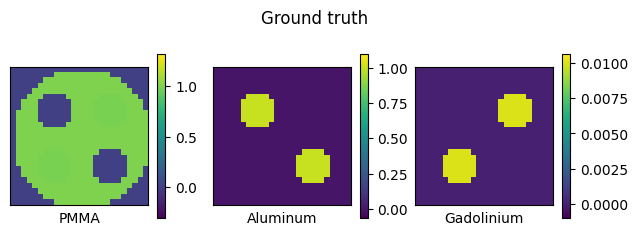}\\\bigskip\bigskip
\includegraphics[width = 0.7\textwidth]{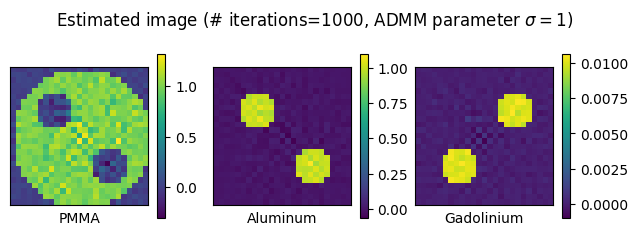}\\\bigskip\bigskip
\includegraphics[width = 0.7\textwidth]{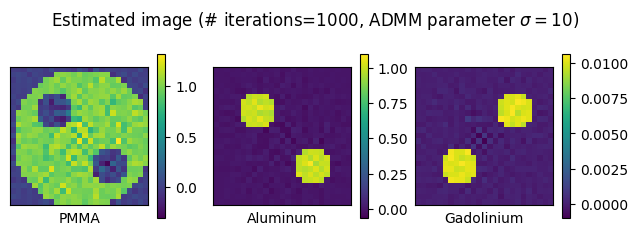}\\\bigskip\bigskip
\includegraphics[width = 0.7\textwidth]{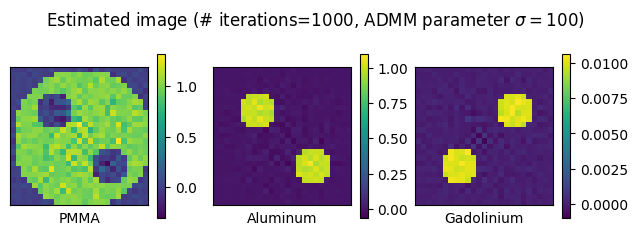}
\caption{The true image in the simulation (top), followed by the reconstructed image (at iteration 1000) 
with each value of the ADMM parameter $\sigma$.
Each row of images
displays the values of $x_{km}$ for each pixel $k$ and each material $m$, for $x = \xh$ (for the ground truth) or $x = x_t$ (for the estimates).}
\label{fig:CT_results_image}
\end{figure}

\begin{figure}\centering
\includegraphics[width = 0.49\textwidth]{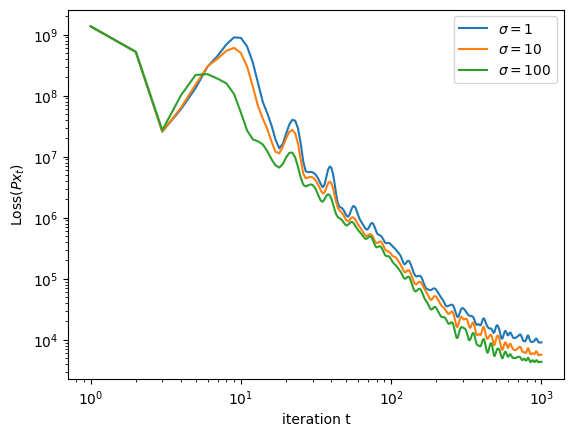}
\includegraphics[width = 0.49\textwidth]{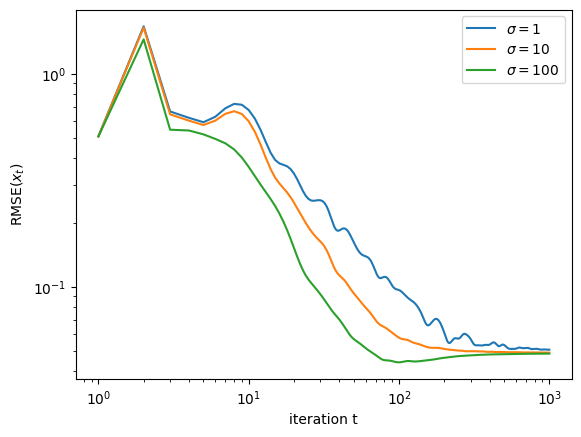}\\\medskip
\includegraphics[width = 0.49\textwidth]{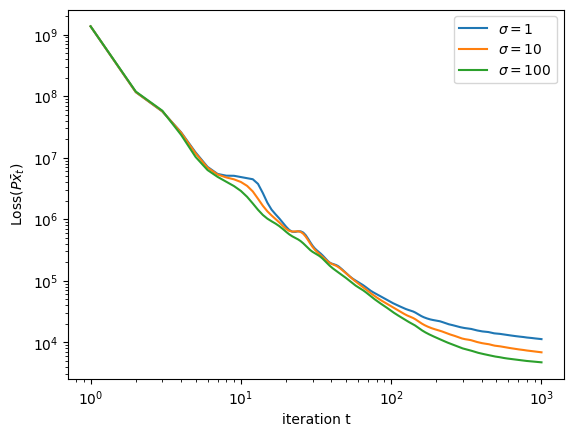}
\includegraphics[width = 0.49\textwidth]{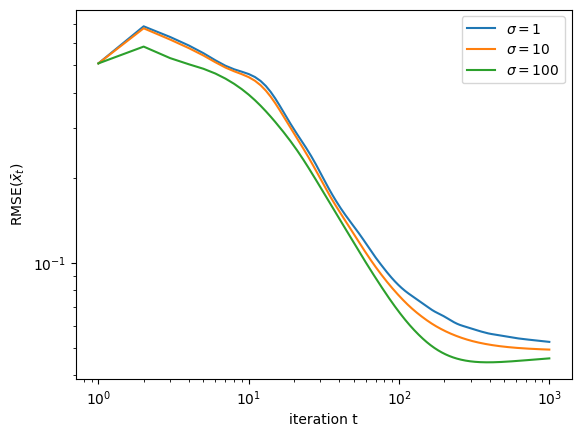}\\
\caption{Convergence results for the CT image reconstruction simulation. 
The top row of the figure shows the value of the objective function $\loss(Px_t)$, and the RMSE $\frac{1}{\sqrt{n_k}}\|x_t - \xh\|_2$, over iteration $t=1,\dots,1000$ of the algorithm,
run with various values of the parameter $\sigma$ as shown. The bottom row shows the same for the running average $\bar{x}_t$ in place of $x_t$.
All axes are on the log scale.}
\label{fig:CT_results_convergence}
\end{figure}

The ground truth, shown in Figure~\ref{fig:CT_results_image}, is a 10cm$\times$10cm two-dimensional image discretized to a $25\times 25$ grid, for a total of $n_k = 25^2=625$ pixels.
The image consists of $n_m=3$ materials---polymethyl methacrylate (PMMA), aluminum, and gadolinium.
As shown in  Figure~\ref{fig:CT_intro},
PMMA  has low attenuation coefficients as it is a plastic, while aluminum, like other metals, has higher attenuation coefficients as it is more difficult for the beam to pass through.
Gadolinium is a contrast material used in clinical CT---its non-monotone attenuation curve
allows for it to be easily identified in the presence of other materials.
The simulated CT scanner has 50 detector cells, and takes images from 50 angles spaced evenly around the unit circle,
for a total of $n_\ell = 50^2 = 2500$ rays along which measurements are taken.
 The beam intensity is set to $10^6$ photons, and there are $n_w = 3$  energy windows, forming a blurry partition of the energy range 
 (see Figure~\ref{fig:CT_intensity}).
 
Figure~\ref{fig:CT_results_image} displays the estimated image (shown at iteration 1000, at each value of the ADMM parameter $\sigma\in\{1,10,100\}$). 
  In Figure~\ref{fig:CT_results_convergence} we show the loss function $\loss(Px_t)$, and the RMSE $\frac{1}{\sqrt{n_k}}\|x_t - \xh\|_2$, at each iteration $t=1,\dots,1000$. As expected, due to the noise in the measurements, the RMSE converges to a small but positive value.
We can see that the algorithm converges steadily towards minimizing the loss  and reducing the RMSE,  and its performance is reasonably stable and robust
across a wide range of values of the tuning parameter $\sigma$.

\paragraph{Extensions} The objective function, and accompanying algorithm, that we have presented here, can easily be modified to incorporate
additional components such as regularization or constraints. 
In particular, total variation regularization can also be incorporated into the framework of
 Algorithm~\ref{algo:main}.\footnote{Details and a demonstration can be found with the code accompanying this paper
  (\url{https://github.com/rinafb/ADMM_CT}), 
  alongside the basic non-regularized simulation setting presented here.  
In addition, this code also shows results from
an  experiment in a noisier setting, with beam intensity set to $10^5$ rather than $10^6$ for a lower signal-to-noise ratio}.
 Another possible modification is adding a preconditioning step to improve the conditioning in the $n_m$-dimensional material 
 space---since the attenuation curves for the three materials are quite similar (see Figure~\ref{fig:CT_intro}),
 adding a preconditioning step can
 improve convergence substantially for the image reconstruction problem (see~\citet{sidky2018three} for more details).
 The algorithm, together with these extensions, has been implemented for large-scale CT data, and has achieved
promising empirical results for both real CT data and simulation
studies, e.g., in \citet{schmidt2022addressing,rizzo2022material,schmidt2023constrained,rizzo2023experimental}.

\begin{figure}\centering
\includegraphics[width = 0.7\textwidth]{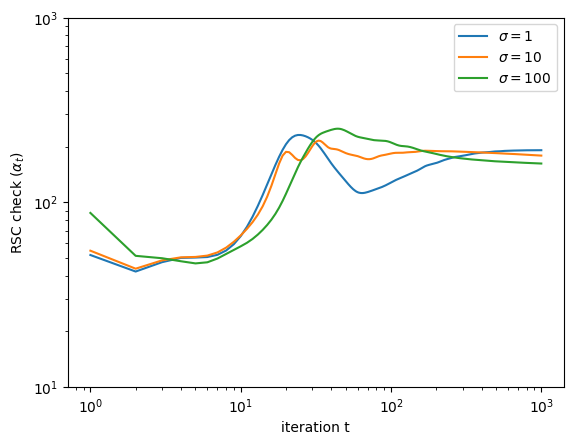}
\caption{The figure shows the value of $\alpha_t$ (defined in~\eqref{eqn:check_rsc_CT_alpha_t}), which empirically 
validates the restricted strong convexity property for the CT imaging example. Both axes are on the log scale.}
\label{fig:CT_results_rsc}
\end{figure}

 \paragraph{Checking assumptions} For the CT imaging example, it is not clear whether it is possible to establish the RSC property~\eqref{eqn:rsc}
theoretically. However, since we are in simulated setting where the target
parameters $(\xh,\yh)$ are known, we can nonetheless validate it empirically. For this example, since $f(x)\equiv 0$, it suffices to check that, for some $\alpha>0$,
\[\inner{y-\yh}{\nabla g(y) - \nabla g(\yh)}\geq \alpha\norm{\vect(y-\yh)}^2_2 - \frac{1}{2}\norm{\vect(Px - y)}^2_\Sigma\]
holds for all $x,y$ (here $\yh=P\xh$). If this is true, then the RSC property holds with $\alpha_f=0$, $\alpha_g=\alpha$, $c_f=c_g=1$,
and $\eps=0$.

However, it is not feasible to verify this over all possible $y$, so we will instead verify that this holds for $y=y_t$ at each iteration $t$ of the algorithm.
 (In fact, examining how the RSC assumption is used in
the proof of Theorem~\ref{thm1}, we see in~\eqref{eqn:use_rsc_in_proof}
that the RSC assumption is only applied at values of $x$ and $y$ appearing along the iterations of the algorithm---specifically,
at points $(x,y)$ of the form $(x_{t+1},y_t)$ at each time $t$. In other words, for the proof of  Theorem~\ref{thm1} to hold for the CT example,
where we have $f(x)\equiv 0$, we only need to check that the inequality above holds at each iteration $y_t$, rather than at all values of $y$.)

To verify this, we calculate
\begin{equation}\label{eqn:check_rsc_CT_alpha_t}\alpha_t := \frac{\inner{ y_t-\yh}{\nabla g( y_t) - \nabla g(\yh)} + \frac{1}{2}\norm{\vect(P x_{t+1} -  y_t)}^2_\Sigma}{\norm{\vect( y_t-\yh)}^2_2},\end{equation}
where $x_{t+1}$ and $y_t$ denote the iterates of the algorithm, while $\yh=P \xh$ is the projection of the true image.
If the RSC property holds as above, then we should see $\alpha_t\geq \alpha$ for all $t$, for some constant $\alpha>0$. Indeed, for the  simulated example, Figure~\ref{fig:CT_results_rsc}
shows that $\alpha_t$ remains bounded away from zero across all iterations of the algorithm. This validates Assumption~\ref{asm:rsc}.

Finally, we verify that approximate first-order optimality~\eqref{eqn:approx_firstorder} holds in this setting.
Choosing $\uh = 0$, we can see that~\eqref{eqn:approx_firstorder} holds 
as long as $\norm{\nabla g(\yh)}_2$ is low. For our simulation, we compare 
$\norm{\nabla g(\yh)}_2$ to $\norm{\nabla g(0)}_2$ (in order for our calculations to be on a meaningful scale),
and we find that
\[\frac{\norm{\nabla g(\yh)}_2}{\norm{\nabla g(0)}_2} =0.000769 ,\]
verifying that approximate first-order optimality holds.

%%%%%%%%%%%%%%%%%%%%%%%%%%%%%%%%%%%%%%%%%%%%%%%%%%%%%%%%%%%%%%%%%%%%%%%%%%%%%

\section{Discussion}\label{sec:discussion} 
The ADMM algorithm has long been known to perform well in a broad range of challenging scenarios, but existing theoretical analyses
are largely restricted to a much more constrained range of settings. Our new theoretical results provide a novel understanding of the performance
of ADMM in the presence of nonsmoothness and nonconvexity in the objective functions, through the lens of a restricted strong convexity property.
A key nonconvex application of this algorithm is the CT image reconstruction problem, where many interesting open questions remain.
In particular, for real CT scanner data, it is important to calibrate the beam intensity and detector sensitivity parameters 
that characterize the performance of the detector. In future work, we aim to extend the ADMM formulation of the image reconstruction problem 
to allow for simultaneous estimation of the calibration parameters (a preliminary study of the simultaneous estimation approach can be found in \citet{ha2018alternating}).  Incorporating more complex aspects of the physical model, such as scatter,
poses an additional challenge that we hope to address in future work to provide a more accurate reconstructed image.

From the theoretical perspective, a key remaining question is whether the RSC property can be further relaxed to allow for convergence
guarantees in an even broader range of settings. On the other hand, the RSC property does not appear to be sufficient to ensure convergence of the
iterates $x_t,y_t$ (rather than the running averages $\bar{x}_t,\bar{y}_t$), as was seen in the quantile regression example. An important open question is whether
a stronger form of the RSC property would allow for convergence guarantees without averaging.
 From the practical side, another important question is the issue of optimization with a stochastic, or mini-batch, approach---analogous to stochastic
gradient descent, the ADMM algorithm can be run using stochastic approximations to gradients at each step (see, e.g., \citet{zhong2014fast}), 
leading to computational speedup, and can be immensely helpful for allowing the method to be applied to large scale applications (including CT imaging, see, e.g., \citet{nien2014fast}).
Another important open question, therefore, is whether the theoretical results of this work for convergence in a nonconvex setting can be extended to the stochastic version of the ADMM algorithm. The empirical performance of the algorithm might also be improved
by incorporating techniques such as adaptive restart \citep{o2015adaptive,kim2018adaptive}, to speed up convergence.

\subsection*{Acknowledgements}
R.F.B. and E.Y.S. were both supported by the National Institutes of Health via grant NIH R01-023968.
R.F.B. was also supported by the National Science Foundation via grant DMS--1654076, and by the Office of Naval Research via grant N00014-20-1-2337. 
E.Y.S. was also supported by National Institutes of Health via grant NIH R01-026282. 
The authors thank Michael Bian for helpful feedback.

\bibliographystyle{plainnat}
\bibliography{bib}

\appendix
\section{Additional details and proofs}

\subsection{A closer look at restricted strong convexity}\label{app:rsc}

To better understand this condition in the setting of the composite optimization problem~\eqref{eqn:optim_intro} studied in this work,
 consider the augmented Lagrangian $\Lcal_{\Sigma}$ defined in~\eqref{eqn:augLagr}.
Since the $x$ and $y$ update steps of ADMM are performing (approximate) alternating minimization on this augmented Lagrangian,
it is intuitive that convexity of the map $(x,y)\mapsto \Lcal_{\Sigma}(x,y,u)$ (at a fixed $u$) is generally needed for convergence to be possible.

On the other hand, if $(x,y)\mapsto \Lcal_{\Sigma/2}(x,y,u)$ is strongly convex (note that we have replaced the penalty matrix $\Sigma$ with 
a smaller penalty, $\Sigma/2$), this is sufficient to ensure the restricted strong convexity condition~\eqref{eqn:rsc} holds (with $\eps=0$)
at any feasible point $(\xh,\yh)$. To see why, 
for any $\xi_x\in\partial f(x)$ and $\zeta_y\in \partial g(y)$, using the fact that $A\xh+B\yh=c$ by feasibility, an elementary calculation 
shows that
\begin{multline}\label{eqn:rsc_vs_Lagrangian}\biginner{\left(\begin{array}{c}x - \xh\\ y - \yh\end{array}\right)}{\left(\begin{array}{c}\xi_x-\xi_{\xh}\\ \zeta_y-\zeta_{\yh}\end{array}\right)} + \frac{1}{2}\norm{Ax + By - c}^2_\Sigma\\
=\biginner{\left(\begin{array}{c}x - \xh\\ y - \yh\end{array}\right)}{\left(\begin{array}{c}\xi_x-\xi_{\xh}\\ \zeta_y-\zeta_{\yh}\end{array}\right)} + \frac{1}{2}\norm{Ax + By - A\xh-B\yh}^2_\Sigma\\
=\biginner{\left(\begin{array}{c}x - \xh\\ y - \yh\end{array}\right)}{\left(\begin{array}{c}\xi_x+\frac{1}{2}A^\top\Sigma (Ax+By)-\xi_{\xh}-\frac{1}{2}A^\top\Sigma (A\xh+B\yh)\\ \zeta_y+\frac{1}{2}B^\top\Sigma (Ax+By)-\zeta_{\yh}-\frac{1}{2}B^\top\Sigma (A\xh+B\yh)\end{array}\right)}.
\end{multline}
We can also calculate 
\[\left(\begin{array}{c}\xi_x+\frac{1}{2}A^\top\Sigma (Ax+By)\\ \zeta_y+\frac{1}{2}B^\top\Sigma (Ax+By)\end{array}\right) \in \partial_{(x,y)}\Lcal_{\Sigma/2}(x,y,\uh)\]
and similarly
\[
\left(\begin{array}{c}\xi_{\xh}+\frac{1}{2}A^\top\Sigma (A\xh+B\yh)\\ \zeta_{\yh}+\frac{1}{2}B^\top\Sigma (A\xh+B\yh)\end{array}\right)\in\partial_{(x,y)}\Lcal_{\Sigma/2}(\xh,\yh,\uh).\]
Therefore, the final expression in~\eqref{eqn:rsc_vs_Lagrangian} will be lower-bounded by strong convexity of $\Lcal_{\Sigma/2}$. Thus, we can interpret the RSC condition~\eqref{eqn:rsc} as only mildly stronger than requiring
strong convexity of the augmented Lagrangian.

\subsection{Completing the proof of Theorem~\ref{thm1}}\label{app:complete_proof_thm1}
To complete the proof of Theorem~\ref{thm1}, we only need to prove that the bound~\eqref{eqn:mainstep} holds under the assumption~\eqref{eqn:stepsize}
on the step size matrices $H_f,H_g$, for any point $(x,y,u)$ with $Ax+By=c$.

By definition of $x_{t+2}$ (i.e., since $x_{t+2}$ is a minimizer of the subproblem that defines its update step),
we must have
\begin{align*}
0&\in\partial f_c(x_{t+2}) + \nabla f_d(x_{t+1}) + A^\top u_{t+1} + A^\top \Sigma (Ax_{t+2}+By_{t+1}-c) + H_f(x_{t+2}-x_{t+1})\\
&=\partial f_c(x_{t+2}) + \nabla f_d(x_{t+1}) + A^\top (2u_{t+1}-u_t)  + A^\top \Sigma A(x_{t+2}-x_{t+1})+ H_f(x_{t+2}-x_{t+1}),\end{align*}
since $u_{t+1}=u_t+\Sigma(Ax_{t+1}+By_{t+1}-c)$.
Since $\partial f(x_{t+2}) = \partial f_c(x_{t+2}) +  \nabla f_d(x_{t+2})$, this implies that there exists some $\xi_{x_{t+2}}\in\partial f(x_{t+2})$ such that
\[\xi_{x_{t+2}}  = \nabla f_d(x_{t+2}) - \nabla f_d(x_{t+1}) - A^\top (2u_{t+1}-u_t) - A^\top \Sigma A(x_{t+2}-x_{t+1}) - H_f(x_{t+2}-x_{t+1})\]
and therefore
\begin{multline*}
\inner{x_{t+2}- x}{\xi_{x_{t+2}} + A^\top u} 
= \inner{x_{t+2}- x}{-A^\top(2u_{t+1} - u_t - u)- A^\top \Sigma A(x_{t+2}-x_{t+1})} \\{}+\inner{x_{t+2}- x}{\nabla f_d(x_{t+2})-\nabla f_d(x_{t+1})  - H_f(x_{t+2}-x_{t+1})} .\end{multline*}
We can similarly calculate
\begin{align*}
0&\in\partial g_c(y_{t+1}) + \nabla g_d(y_t) + B^\top u_t + B^\top \Sigma (Ax_{t+1}+By_{t+1}-c) + H_g(y_{t+1}-y_t)\\
&=\partial g_c(y_{t+1})   + \nabla g_d(y_t) + B^\top u_{t+1} + H_g(y_{t+1}-y_t),\end{align*}
and so there exists some $\zeta_{y_{t+1}}\in\partial g(y_{t+1})$ satisfying
\begin{multline*}
\inner{y_{t+1}- y}{\zeta_{y_{t+1}} + B^\top u} = \inner{y_{t+1}- y}{  -B^\top ( u_{t+1}-u )} \\
{}+\inner{y_{t+1}- y}{\nabla g_d(y_{t+1})-\nabla g_d(y_t)  - H_g(y_{t+1}-y_t)}.\end{multline*}
We can further calculate
\begin{align*}
&\inner{x_{t+2}- x}{-A^\top(2u_{t+1} - u_t - u)- A^\top \Sigma A(x_{t+2}-x_{t+1})} \\
&=\left(\begin{array}{c} x_{t+2}-x\\ u_{t+1}-u\end{array}\right)^\top \left(\begin{array}{cc} A^\top \Sigma A & A^\top \\  A & \Sigma^{-1}\end{array}\right) \left(\begin{array}{c} x_{t+1} - x_{t+2} \\ u_t - u_{t+1}\end{array}\right) - \inner{\Sigma^{-1}(u_t-u_{t+1})+A(x_{t+1}-x)}{u_{t+1}-u}\\
&=\left(\begin{array}{c} x_{t+2}-x\\ u_{t+1}-u\end{array}\right)^\top \left(\begin{array}{cc} A^\top \Sigma A& A^\top \\  A & \Sigma^{-1}\end{array}\right) \left(\begin{array}{c} x_{t+1} - x_{t+2} \\ u_t - u_{t+1}\end{array}\right) +\inner{B(y_{t+1}-y)}{u_{t+1}-u}.
 \end{align*}
 
Combining our calculations so far, we have
\begin{multline}\label{eqn:rearrange_step1}
 \biginner{\left(\begin{array}{c}x_{t+2}-x\\ y_{t+1}-y\end{array}\right)}{\left(\begin{array}{c}\xi_{x_{t+2}} +A^\top u\\\zeta_{y_{t+1}} +B^\top u\end{array}\right)}
= (z_{t+1}-z)^\top M(z_t-z_{t+1}) \\ {}+ \inner{x_{t+2}- x}{\nabla f_d(x_{t+2})-\nabla f_d(x_{t+1})} + \inner{y_{t+1}- y}{\nabla g_d(y_{t+1})-\nabla g_d(y_t)} ,\end{multline}
where we define
$z = (x,y,u)$ and $z_t = (x_{t+1},y_t,u_t)$ for each $t$, and let
\[M =  \left(\begin{array}{ccc}H_f+A^\top\Sigma A& 0 & A^\top \\ 0& H_g & 0 \\ A &   0 & \Sigma^{-1}\end{array}\right)\succeq 0.\]
Next, defining $\norm{v}_M = \sqrt{v^\top Mv}$, we can use a telescoping sum to calculate
\begin{multline*}
\sum_{t=0}^{T-1}(z_{t+1}-z)^\top M(z_t-z_{t+1}) 
=\sum_{t=0}^{T-1} \left(\frac{1}{2}\norm{z - z_t}^2_M-\frac{1}{2}\norm{z- z_{t+1}}^2_M-\frac{1}{2}\norm{z_t - z_{t+1}}^2_M\right)\\
=\frac{1}{2}\norm{z - z_0}^2_M - \frac{1}{2}\norm{z - z_T}^2_M - \frac{1}{2}\sum_{t=0}^{T-1} \norm{z_t - z_{t+1}}^2_M.
\end{multline*}
Furthermore,
\begin{multline*}\norm{z_t - z_{t+1}}^2_M - \norm{x_{t+1}-x_{t+2}}^2_{H_f} - \norm{y_t-y_{t+1}}^2_{H_g}\\
=\left(\begin{array}{c} x_{t+1}-x_{t+2}\\ u_t - u_{t+1}\end{array}\right)^\top \left(\begin{array}{cc} A^\top \Sigma A & A^\top \\  A & \Sigma^{-1}\end{array}\right) \left(\begin{array}{c} x_{t+1} - x_{t+2}\\ u_t - u_{t+1}\end{array}\right)\\ = \norm{\Sigma^{-1}(u_t-u_{t+1})+A(x_{t+1}-x_{t+2})}^2_\Sigma = \norm{Ax_{t+2}+By_{t+1}-c}^2_\Sigma,\end{multline*}
where the last step plugs in the update step for $u_{t+1}$. 
Combining these calculations with~\eqref{eqn:rearrange_step1}, we obtain
\begin{multline}\label{eqn:rearrange_step2}
\sum_{t=0}^{T-1} \biginner{\left(\begin{array}{c}x_{t+2}-x\\ y_{t+1}-y\end{array}\right)}{\left(\begin{array}{c}\xi_{x_{t+2}} +A^\top u\\\zeta_{y_{t+1}} +B^\top u\end{array}\right)} \\
= \frac{1}{2}\norm{z - z_0}^2_M - \frac{1}{2}\norm{z - z_T}^2_M -\frac{1}{2}\sum_{t=0}^{T-1} \norm{Ax_{t+2}+By_{t+1}-c}^2_\Sigma\\
{}+\sum_{t=0}^{T-1} \left[\inner{x_{t+2}- x}{\nabla f_d(x_{t+2})-\nabla f_d(x_{t+1})} - \frac{1}{2} \norm{x_{t+1}-x_{t+2}}^2_{H_f} \right] \\{}+ \sum_{t=0}^{T-1}\left[\inner{y_{t+1}- y}{\nabla g_d(y_{t+1})-\nabla g_d(y_t)} - \frac{1}{2}\norm{y_t-y_{t+1}}^2_{H_g}\right].\end{multline}

Now, since $H_f\succeq \nabla^2 f_d(x)$ by the assumption~\eqref{eqn:stepsize},  we can write
\[f_d(x_{t+2}) \leq f_d(x_{t+1}) + \inner{x_{t+2}-x_{t+1}}{\nabla f_d(x_{t+1})} + \frac{1}{2}\norm{x_{t+1}-x_{t+2}}^2_{H_f}\]
for each $t$. Rearranging terms and taking a telescoping sum, this means that
\begin{multline*}
 \sum_{t=0}^{T-1}\left[\inner{x_{t+2}- x}{\nabla f_d(x_{t+2})-\nabla f_d(x_{t+1})} - \frac{1}{2} \norm{x_{t+1}-x_{t+2}}^2_{H_f} \right]
\\
\leq  f_d(x_1) - f_d(x_{T+1}) +  \inner{x-x_1}{\nabla f_d(x_1)} - \inner{x-x_{T+1}}{\nabla f_d(x_{T+1})} .\end{multline*}
Again applying $H_f\succeq \nabla^2 f_d(x)$, we also have
\[f_d(x)\leq f_d(x_{T+1}) + \inner{x-x_{T+1}}{\nabla f_d(x_{T+1})}+\frac{1}{2}\norm{x-x_{T+1}}^2_{H_f}\]
and
\[f_d(x_1)\leq f_d(x) + \inner{x_1-x}{\nabla f_d(x)}+\frac{1}{2}\norm{x-x_1}^2_{H_f},\]
which combined with the above yields
\begin{multline*}
 \sum_{t=0}^{T-1}\left[\inner{x_{t+2}- x}{\nabla f_d(x_{t+2})-\nabla f_d(x_{t+1})} - \frac{1}{2} \norm{x_{t+1}-x_{t+2}}^2_{H_f} \right]
\\
\leq  - \inner{x-x_1}{\nabla f_d(x) - \nabla f_d(x_1)} + \frac{1}{2}\norm{x-x_{T+1}}^2_{H_f}+\frac{1}{2}\norm{x-x_1}^2_{H_f} .\end{multline*}
Performing an identical calculation for the $y$ terms,
and combining these calculations with~\eqref{eqn:rearrange_step2} along with the fact that $\norm{z - z_T}^2_M\geq\norm{x-x_{T+1}}^2_{H_f}+\norm{y-y_T}^2_{H_g}$, we obtain
\begin{multline*}\sum_{t=0}^{T-1} \biginner{\left(\begin{array}{c}x_{t+2}-x\\ y_{t+1}-y\end{array}\right)}{\left(\begin{array}{c}\xi_{x_{t+2}} +A^\top u\\\zeta_{y_{t+1}} +B^\top u\end{array}\right)}\\ \leq  C_1(x,y,u;x_0,y_0,u_0)- \frac{1}{2}\sum_{t=0}^{T-1} \norm{Ax_{t+2}+By_{t+1}-c}^2_\Sigma,\end{multline*}
where we define
\begin{multline}\label{eqn:mainstep_define_C}
 C_1(x,y,u;x_0,y_0,u_0)
= \frac{1}{2}\norm{z - z_0}^2_M 
- \inner{x-x_1}{\nabla f_d(x) - \nabla f_d(x_1)} +\frac{1}{2}\norm{x-x_1}^2_{H_f} \\
 {}  - \inner{y-y_0}{\nabla g_d(y) - \nabla g_d(y_0)} +\frac{1}{2}\norm{y-y_0}^2_{H_g} .\end{multline}
 (Note that $x_1$ is a deterministic function of $(x_0,y_0,u_0)$, and therefore $C_1$ can depend implicitly on $x_1$.)
This proves the desired bound~\eqref{eqn:mainstep}.

\subsection{Details for implementing ADMM for the CT application}\label{app:CT}

To run Algorithm~\ref{algo:main} for the CT image reconstruction problem~\eqref{eqn:CT_opt_xy}, plugging in our choices of parameters $H_f,H_g,\Sigma$
and the values of $A,B,c$ and $f(x)\equiv 0$, our update steps can be calculated as follows. Note that in our notation below,
the $x,y,u$ variables are all treated as matrices, with $n_k\times n_m$ dimensional $x$ variables
and with $n_\ell\times n_m$ dimensional $y$ and $u$ variables.
\begin{itemize}
\item The $x$ update step is given by
\[x_{t+1} =x_t + Q_f^{-1}P^\top \big(\tilde\Sigma(y_t - Px_t) - u_t\big).\]
 Since $Q_f$ and $\tilde\Sigma$ are diagonal while $P$ is sparse, this requires only inexpensive
 matrix-vector calculations.
 \item The $y$ update step is given by solving the minimization problem
 \[y_{t+1} = \argmin_y\left\{g_c(y) + \inner{y}{\nabla g_d(y_t) - ( u_t + \tilde\Sigma Px_{t+1}  )}+ \frac{1}{2}\vect(y)^\top (\tilde\Sigma\otimes\ident_{n_m}) \vect(y) \right\}.\]
We recall from the definition of $g_c$~\eqref{eqn:CT_define_g_c} that this function separates over the $n_\ell$ many rays---that is, we
can write $g_c(y) = \sum_\ell g_{c,\ell}(y_\ell)$, where $y_\ell\in\R^{n_m}$ is the portion of $y$ corresponding to the $\ell$-th ray,
and where
\[g_{c,\ell}(y_\ell) =  \sum_{w} \sum_i S_{w\ell i} \qexp\left\{ - \sum_m\mu_{mi}y_{\ell m}\right\}.\]
Therefore, equivalently, the $y$ update step is given by solving
 \[(y_{t+1})_\ell = \argmin_{y_\ell\in\R^{n_m}}\left\{g_{c,\ell}(y) + \inner{y_\ell}{(\nabla g_d(y_t))_\ell - (u_t)_\ell - \tilde\Sigma_{\ell\ell}( Px_{t+1}  )_\ell}+ \frac{\tilde\Sigma_{\ell\ell}}{2}\norm{y_\ell}^2_2\right\}\]
for each $\ell=1,\dots,n_\ell$. Since we typically work with a small number of materials (e.g., $n_m=3$ or $n_m=5$), solving each one of 
these convex minimization problems
 is computationally very inexpensive.
We will use  the Newton--Raphson method 
to solve the minimization subproblem approximately, in parallel for each $\ell$:
setting $y_{t+1}^{(0)} = y_t$, we define
\begin{multline*}
(y_{t+1}^{(i+1)})_{\ell} = (y_{t+1}^{(i)})_{\ell} -(\nabla^2 g_{c,\ell}((y_{t+1}^{(i)})_\ell) + \tilde\Sigma_{\ell\ell}\ident_{n_m})^{-1} \cdot{}\\\Big(\nabla g_{c,\ell}((y_{t+1}^{(i)})_\ell) + (\nabla g_d(y_t))_\ell + \tilde\Sigma_{\ell\ell}(y_{t+1}^{(i)} - Px_{t+1})_\ell - (u_t)_\ell\Big),\end{multline*}
for each $i=0,1,2,\dots,N-1$, and then set $y_{t+1} = y_{t+1}^{(N)}$.  In our implementation,
at each iteration $t$ we run $N=10$ steps of the Newton--Raphson method  to compute the $y$ update, which is sufficient to obtain a near-exact solution.

\item The $u$ update step is given by
\[u_{t+1} =  u_t + \tilde\Sigma (P x_{t+1} - y_{t+1}).\]
 Since $\tilde\Sigma$ is diagonal while $P$ is sparse, this again requires only inexpensive
 matrix-vector calculations.
 \end{itemize}

\subsection{Details for implementing ADMM for the sparse quantile regression example}\label{app:sparseQR}
We now compute  the steps of Algorithm~\ref{algo:main} for the sparse quantile regression example, i.e., for the problem
of minimizing~\eqref{eqn:sparseQR}.
Plugging in our choices of the parameters $H_f,H_g,\Sigma$ and of $A,B,c$, the steps of Algorithm~\ref{algo:main} are given by
\begin{align*}
x_{t+1} &= \argmin_{x\in\R^d} \left\{f_c(x) +  \inner{x}{ \nabla f_d(x_t) +\Phi^\top u_t}+ \frac{\sigma}{2}\norm{\Phi x - y_t}^2_2 + \frac{\sigma}{2}\norm{x - x_t}^2_{\gamma\ident_d - \Phi^\top \Phi}\right\}, \\
y_{t+1} &= \argmin_{y\in\R^n}\left\{g_c(y) +  \inner{y}{ \nabla g_d(y_t)-u_t } + \frac{\sigma}{2}\norm{\Phi x_{t+1} - y}^2_2\right\},\\
u_{t+1} &= u_t + \sigma(\Phi x_{t+1}-y_{t+1}).
\end{align*}
Now we compute the $x$ and $y$ update steps explicitly. First, for $x$, recall that $f_c(x)=\lambda\norm{x}_1+ \delta_{\norm{x}_2\leq R}$ and
\[f_d(x) = \lambda\sum_{j=1}^d \left(\beta\log(1+ |x_j|/\beta) - |x_j|\right).\]
We can calculate the gradient as
\[\big[\nabla f_d(x)\big]_j  = -\frac{\lambda x_j}{\beta+|x_j|}.\] 
Therefore,
\begin{multline*}
f_c(x) +  \inner{x}{ \nabla f_d(x_t) +\Phi^\top u_t}+ \frac{\sigma}{2}\norm{\Phi x - y_t}^2_2 + \frac{\sigma}{2}\norm{x - x_t}^2_{\gamma\ident_d - \Phi^\top \Phi} \\
=\sum_{j=1}^d \left( \frac{\sigma\gamma}{2}x_j^2 - x_j \cdot\left[\frac{\lambda (x_t)_j}{\beta+|(x_t)_j|} + \sigma (\gamma x_t - \Phi^\top (\Phi x_t-y_t +u_t/\sigma))_j\right]\right) 
+\lambda \norm{x}_1+ \delta_{\norm{x}_2\leq R}\\
=\frac{\sigma\gamma}{2}\sum_{j=1}^d \left( x_j^2 -2  x_j \cdot (\tilde{x}_{t+1})_j\right) 
+\lambda \norm{x}_1+ \delta_{\norm{x}_2\leq R},
\end{multline*}
where we define a vector $\tilde{x}_{t+1}$ with entries
\[(\tilde{x}_{t+1})_j =(x_t)_j - \frac{\big(\Phi^\top (\Phi x_t-y_t +u_t/\sigma)\big)_j}{\gamma}+ \frac{\lambda}{\sigma\gamma}\cdot\frac{(x_t)_j}{\beta+|(x_t)_j|}.\]
Then we can verify that the objective function above is minimized by defining
\[x_{t+1} =  \textnormal{SoftThresh}_{\frac{\lambda}{\sigma\gamma}}\left(\tilde{x}_{t+1}\right) \cdot\min\left\{1 , \frac{R}{\norm{\textnormal{SoftThresh}_{\frac{\lambda}{\sigma\gamma}}\left(\tilde{x}_{t+1}\right)}_2}\right\},\]
where the soft thresholding function, $\textnormal{SoftThresh}_\lambda:\R^d\rightarrow\R^d$, is defined elementwise as 
\[\big[\textnormal{SoftThresh}_\lambda(x) \big]_j= \begin{cases} x_j - \lambda, & \textnormal{ if }x_j > \lambda, \\ 0 , & \textnormal{ if } |x_j|\leq \lambda, \\ x_j + \lambda, & \textnormal{ if } x_j < -\lambda.\end{cases}\]

Next, for the $y$ update step, recall $g_d(y)\equiv 0$ and
\[g_c(y) = \frac{1}{n} \sum_{i=1}^n q\max\{w_i - y_i,0\} + (1-q)\max\{y_i-w_i,0\}.\]
Then the optimization problem for the $y$ update step  separates over the $n$ entries of $y$:
\begin{multline*}
g_c(y) +  \inner{y}{ \nabla g_d(y_t)-u_t } + \frac{\sigma}{2}\norm{\Phi x_{t+1} - y}^2_2\\
=\sum_{i=1}^n \left(\frac{1}{n}\left[q\max\{w_i - y_i,0\} + (1-q)\max\{y_i-w_i,0\}\right] + \frac{\sigma}{2}y_i^2 -  y_i \cdot \left(\sigma (\Phi x_{t+1})_i + (u_t)_i\right)\right).\end{multline*}
This is minimized by setting $y_{t+1}$ to have entries
\[(y_{t+1})_i = \begin{cases}
(\Phi x_{t+1})_i + \frac{(u_t)_i}{\sigma} + \frac{q}{n\sigma}, & \textnormal{ if } (\Phi x_{t+1})_i + \frac{(u_t)_i}{\sigma}  +  \frac{q}{n\sigma} < w_i,\\
(\Phi x_{t+1})_i + \frac{(u_t)_i}{\sigma}  -\frac{1-q}{n\sigma}, & \textnormal{ if } (\Phi x_{t+1})_i + \frac{(u_t)_i}{\sigma}  - \frac{1-q}{n\sigma} > w_i,\\
w_i, & \textnormal{ if } (\Phi x_{t+1})_i + \frac{(u_t)_i}{\sigma}  -  \frac{1-q}{n\sigma}\leq w_i \leq (\Phi x_{t+1})_i + \frac{(u_t)_i}{\sigma}  + \frac{q}{n\sigma} .
\end{cases}\]

\subsection{Proof of Proposition~\ref{prop:quantile_rsc} (verifying assumptions for the sparse quantile regression example)}\label{app:sparseQR_proof}
To prove the result, we need to check that, with  probability at least $1-(nd)^{-1}$, the RSC bound~\eqref{eqn:rsc} 
and the approximate first-order optimality condition~\eqref{eqn:approx_firstorder} both  at the point $(\xh,\yh,\uh)$,
with parameters defined as in the statement of the proposition. Concretely, 
let $\uh\in\R^n$ have entries 
\[\uh_i = \frac{1}{n}(-q\cdot \One{z_i>0} + (1-q)\cdot\One{z_i<0}).\]
Then we can verify $\zeta_{\yh} = \uh \in\partial g(\yh)$. Define also $\xi_{\xh}$ to have entries
\[(\xi_{\xh})_j = \begin{cases}\frac{\lambda \beta \sign(\xh_j)}{\beta +|\xh_j|} , &\xh_j\neq 0,\\ (-\Phi^\top \uh)_j , & \xh_j=0.\end{cases}\]

We will show that, with the desired probability,
\begin{multline}\label{eqn:quantile_rsc_to_show}\biginner{\left(\begin{array}{c}x - \xh\\ y - \yh\end{array}\right)}{\left(\begin{array}{c} \xi_x -\xi_{\xh} \\ \zeta_y-\zeta_{\yh}\end{array}\right)} \\
\geq C_3\min\left\{\norm{x-\xh}^2_2,\norm{x-\xh}_2\right\} - \frac{\sigma}{2}\norm{y-\Phi x}^2_2 - C_4\max\{1,\sigma^{-1}\}\cdot\frac{s_*\log(nd)}{n}\end{multline}
for all $x\in\dom(f)$ $y\in\dom(g)$, $\xi_x\in\partial f(x)$, $\zeta_y\in\partial g(y)$, and that
\begin{equation}\label{eqn:quantile_firstorder_to_show}\xi_{\xh}\in\partial f(\xh)\text{ and }\norm{-\Phi^\top \uh - \xi_{\xh}}_2\leq \min\left\{\frac{C_3}{2},\sqrt{C_3 C_5}\sqrt{\frac{s_*\log(nd)}{n}}\right\},\end{equation}
where $C_3,C_4,C_5>0$ are constants that depend only on $c_z,t_z,a_\phi,b_\phi,B_\phi$ and on $C_\lambda$.
These bounds are sufficient to verify Assumptions~\ref{asm:rsc} and~\ref{asm:approx_firstorder}, as desired.

\subsubsection{Verifying approximate first-order optimality}
First we check that $\xi_{\xh}\in\partial f(\xh)$.
Recall that we can write
$f(x) = f_c(x) + f_d(x)$ where, for any $x\in\dom(f)$ (i.e., $\norm{x}_2\leq R$), we have \[f(x)=   \lambda\sum_{j=1}^d \beta\log(1 + |x_j|/\beta).\]
Now fix any $x\in\dom(f)$. Then we can calculate that a subgradient $\xi_x\in\partial f(x)$ must have entries satisfying
\begin{equation}\label{eqn:calculate_subgradient_QR}\begin{cases}(\xi_x)_j=\frac{\lambda\beta\sign(x_j)}{\beta+|x_j|}, & x_j \neq 0,\\
(\xi_x)_j\in [-\lambda,\lambda], & x_j = 0.\end{cases}\end{equation}
From this calculation, we can see that to verify $\xi_{\xh}\in\partial f(\xh)$, we only need to check that $|(\xi_{\xh})_j|\leq \lambda$ for all $j$ with $\xh_j=0$.
Since $\norm{\Phi}_{\infty}\leq B_\phi$ with probability 1, while $\uh$ is a $\frac{1}{n}$-bounded zero-mean vector, 
Hoeffding's inequality shows that
\begin{equation}\label{eqn:QR_Hoeffding}\PP{ \norm{\Phi^\top \uh}_{\infty} \leq 2B_\phi \sqrt{\frac{\log(nd)}{n}}} \geq 1-(2nd)^{-1}.\end{equation}
From this point on, we will assume that this event holds.
Since $\lambda  = C_\lambda  \sqrt{\frac{\log(nd)}{n}} \geq C_1  \sqrt{\frac{\log(nd)}{n}} \geq 2B_\phi \sqrt{\frac{\log(nd)}{n}}$ (as long as we take $C_1 \geq 2B_\phi$, as we will do below),
this verifies that $|(\xi_{\xh})_j|\leq \lambda$ for $j$ such that $\xh_j=0$, and thus 
$\xi_{\xh}\in\partial f(\xh)$, as desired.

Next we check that~\eqref{eqn:quantile_firstorder_to_show} holds, to complete our verification
of the approximate first-order optimality assumption.
Writing $S_*\subseteq\{1,\dots,d\}$ to denote the support of $\xh$, we have
\[\norm{(\xi_{\xh})_{S_*}}_2 = \sqrt{\sum_{j:\xh_j\neq 0}\left(\frac{\lambda \beta \sign(\xh_j)}{\beta +|\xh_j|} \right)^2 }\leq \sqrt{\sum_{j:\xh_j\neq 0}\lambda^2 }\leq \sqrt{s_*}\lambda,\]
and also,
\[ \norm{(\Phi^\top \uh)_{S_*}}_2  \leq \sqrt{s_*}\norm{(\Phi^\top \uh)_{S_*}}_{\infty} \leq \sqrt{s_*}\cdot  2B_\phi \sqrt{\frac{\log(nd)}{n}}.\]
Then
\begin{multline*}\norm{-\Phi^\top \uh - \xi_{\xh}}_2 \leq \norm{(-\Phi^\top \uh - \xi_{\xh})_{S_*}}_2  + \norm{(-\Phi^\top \uh - \xi_{\xh})_{S^c_*}}_2\\
 \leq \norm{(\Phi^\top \uh)_{S_*}}_2 + \norm{(\xi_{\xh})_{S_*}}_2  + \norm{(-\Phi^\top \uh - \xi_{\xh})_{S^c_*}}_2.\end{multline*}
Since $ \norm{(-\Phi^\top \uh - \xi_{\xh})_{S^c_*}}_2=0$ by definition, this establishes that
\[\norm{-\Phi^\top \uh - \xi_{\xh}}_2 \leq \sqrt{s_*}\left(\lambda +  2B_\phi \sqrt{\frac{\log(nd)}{n}}\right) = (C_\lambda+2B_\phi)\sqrt{\frac{s_*\log(nd)}{n}}.\]
Recall $\frac{s_*\log(nd)}{n}\leq C_0$ by assumption, and furthermore,
\begin{equation}\label{eqn:bound_lambda_root_s}\lambda\sqrt{s_*} = C_\lambda \sqrt{\frac{\log(nd)}{n}} \cdot s_*^{1/2} \leq C_1\sqrt{\frac{C_0 \frac{n}{\log(nd)}}{s_*}}  \sqrt{\frac{\log(nd)}{n}} \cdot s_*^{1/2}
  = C_1\sqrt{C_0}.\end{equation}
  We therefore have
\[\norm{-\Phi^\top \uh - \xi_{\xh}}_2 \leq\min\left\{\sqrt{C_0}(C_1+2B_\phi), (C_\lambda+2B_\phi)\sqrt{\frac{s_*\log(nd)}{n}}\right\}.\]
Finally, choosing the constants as $C_3 = 2\sqrt{C_0}(C_1 + 2B_\phi)$ and $C_5 =  (C_\lambda + 2B_\phi)^2/C_3$,
we have proved~\eqref{eqn:quantile_firstorder_to_show}.
\subsubsection{Verifying restricted strong convexity}
Next we will verify that~\eqref{eqn:quantile_rsc_to_show} holds, to validate the restricted
strong convexity property.

\paragraph{Bounding the $x$ term}
Recall our earlier calculation~\eqref{eqn:calculate_subgradient_QR} of the subgradient $\partial f(x)$.
Writing $S_*\subseteq\{1,\dots,d\}$ to denote the support of $\xh$ as before, for each $j\in S_*^c$ we have
\[(x-\xh)_j \cdot (\xi_x)_j = x_j\cdot(\xi_x)_j = \frac{\lambda \beta |x_j|}{\beta+|x_j|} \geq \lambda|x_j|-\lambda\beta^{-1}x_j^2\]
if $x_j\neq 0$, or if $x_j=0$ then $(x-\xh)_j \cdot (\xi_x)_j = 0 = \lambda|x_j|-\lambda\beta^{-1}x_j^2$ holds trivially.
Thus
\[\inner{(x-\xh)_{S_*^c}}{(\xi_x)_{S_*^c}}  \geq  \lambda\norm{x_{S_*^c}}_1- \lambda\beta^{-1}\norm{x_{S_*^c}}^2_2= \lambda\norm{(x-\xh)_{S_*^c}}_1- \lambda\beta^{-1}\norm{(x-\xh)_{S_*^c}}^2_2.\]
Next, since $(\xi_{\xh})_{S_*^c} = (-\Phi^\top\uh)_{S_*^c}$ and we know that  $ \norm{\Phi^\top \uh}_{\infty} \leq 2B_\phi \sqrt{\frac{\log(nd)}{n}}$ by~\eqref{eqn:QR_Hoeffding}, we have
\[\inner{(x-\xh)_{S_*^c}}{(\xi_x-\xi_{\xh})_{S_*^c}}  \geq \left(  \lambda-2B_\phi \sqrt{\frac{\log(nd)}{n}}\right)\norm{x_{S_*^c}}_1- \lambda\beta^{-1}\norm{x_{S_*^c}}^2_2 .\]

Next,  the function $t\mapsto \beta\log(1+|t|/\beta)$ can be decomposed as
\[ \beta\log(1+|t|/\beta) = |t| + \left( \beta\log(1+|t|/\beta) - |t|\right),\]
where the first term is convex while the second term is concave and twice differentiable with second derivative $\geq -\beta^{-1}$, 
which proves that
\[\inner{(x-\xh)_{S_*}}{(\xi_x-\xi_{\xh})_{S_*}} \geq -\lambda\beta^{-1}\norm{(x-\xh)_{S_*}}^2_2.\] 

Putting all our calculations together, we have established that
\begin{multline*}\inner{x-\xh}{\xi_x-\xi_{\xh}} \geq \left(  \lambda-2B_\phi \sqrt{\frac{\log(nd)}{n}}\right)\norm{x_{S_*^c}}_1- \lambda\beta^{-1}\norm{x-\xh}^2_2\\
 \geq  \left(\lambda - 2B_\phi \sqrt{\frac{\log(nd)}{n}} \right)\norm{x-\xh}_1 - \lambda\beta^{-1}\norm{x-\xh}^2_2 -\lambda \norm{(x-\xh)_{S_*}}_1\\
  \geq  (C_1 - 2B_\phi) \sqrt{\frac{\log(nd)}{n}} \norm{x-\xh}_1 - \lambda\beta^{-1}\norm{x-\xh}^2_2 -\lambda s_*^{1/2} \norm{x-\xh}_2,\end{multline*}
  where the last step holds since  $|S_*|\leq s_*$, and by definition of $\lambda$.
 Finally, if $\norm{x-\xh}_2\leq 1$, then we have
 \begin{multline*}
 \lambda\beta^{-1}\norm{x-\xh}^2_2 +\lambda s_*^{1/2} \norm{x-\xh}_2
 \leq 
C_2^{-1}\norm{x-\xh}^2_2 +C_\lambda \sqrt{\frac{s_*\log(nd)}{n}}\norm{x-\xh}_2\\
  \leq 2C_2^{-1}\norm{x-\xh}^2_2 + \frac{C_2 C_\lambda^2}{4}\cdot\frac{s_*\log(nd)}{n},
  \end{multline*}
 by our bound on $\beta$ along with the fact that $ab\leq ca^2/2 + b^2/2c$ for all $a,b,c>0$. If instead $\norm{x-\xh}_2> 1$,
 then $\norm{x-\xh}_2\leq 2R$ since $x,\xh\in\dom(f)$, and so
 \begin{multline*}
 \lambda\beta^{-1}\norm{x-\xh}^2_2 +\lambda s_*^{1/2} \norm{x-\xh}_2
 \leq \left(  2\lambda\beta^{-1}R + \lambda s_*^{1/2}\right)\norm{x-\xh}_2\\
 \leq \left(2C_2^{-1} + C_1\sqrt{C_0}\right)\norm{x-\xh}_2
  \end{multline*}
 since $2\lambda\beta^{-1}R\leq 2C_2^{-1}$ by our bound on $\beta$,
  and since $\lambda s_*^{1/2}\leq C_1\sqrt{C_0}$ as calculated in~\eqref{eqn:bound_lambda_root_s} above.
 Therefore, combining everything,
 \begin{multline*}\inner{x-\xh}{\xi_x-\xi_{\xh}} \geq(C_1 - 2B_\phi) \sqrt{\frac{\log(nd)}{n}} \cdot \norm{x-\xh}_1\\
 - \left(2C_2^{-1} + C_1\sqrt{C_0} \right) \cdot\min\{\norm{x-\xh}^2_2,\norm{x-\xh}_2\}
 -  \frac{C_2 C_\lambda^2}{4}\cdot\frac{s_*\log(nd)}{n}.\end{multline*}
 
\paragraph{Bounding the $y$ term}
First, we compute the subgradient of $t\mapsto \ell_q(s-t)$:
\[\partial_t \ell_q(s-t) = \begin{cases} \{-q\}, &  t < s,\\
[-q,1-q], & t=s,\\
\{1-q\}, & t>s.\end{cases}\]
Therefore any $\zeta_y \in \partial g(y)$ must have entries satisfying
\[ \begin{cases}  n(\zeta_y)_i =  -q, & y_i < \yh_i + z_i,\\
 n(\zeta_y)_i \in [-q,1-q], & y_i = \yh_i + z_i,\\
 n(\zeta_y)_i =1-q, & y_i > \yh_i + z_i.\end{cases}\]
By definition of $\zeta_{\yh}$ from above, we can therefore calculate
\[ \begin{cases} n\big((\zeta_y)_i - (\zeta_{\yh})_i \big)= 
0, & \text{ if $z_i>0$ and $z_i>y_i - \yh_i$, or $z_i<0$ and $z_i<y_i - \yh_i$},\\
n\big((\zeta_y)_i - (\zeta_{\yh})_i\big) \in[0,1], & \text{ if $z_i>0$ and $z_i=y_i - \yh_i$},\\
n\big((\zeta_y)_i - (\zeta_{\yh})_i\big)  \in[-1,0], & \text{ if $z_i<0$ and $z_i=y_i - \yh_i$},\\
n\big((\zeta_y)_i - (\zeta_{\yh})_i\big) = -1, & \text{ if $z_i <0$ and $z_i > y_i-\yh_i$},\\
n\big((\zeta_y)_i - (\zeta_{\yh})_i\big) = 1,  & \text{ if $z_i >0$ and $z_i < y_i-\yh_i$}.
\end{cases}\]
(Note that $z_i\neq 0$ almost surely, so we can ignore the case $z_i=0$.)
We can therefore calculate
\begin{multline*}\inner{y-\yh }{\zeta_y - \zeta_{\yh}} \geq \frac{1}{n}\cdot \overbrace{\sum_i (y_i - \yh_i) \cdot \One{ y_i-\yh_i>z_i>0}}^{\textnormal{Term 1}} \\{}+\frac{1}{n}\cdot \underbrace{\sum_i (\yh_i - y_i) \cdot \One{\yh_i-y_i>-z_i>0}}_{\textnormal{Term 2}} .\end{multline*}

Writing $(t)_+ = \max\{t,0\}$ for any $t\in\R$, we then have
\begin{align*}
&\textnormal{Term 1}
=\sum_i (y_i - \yh_i) \cdot \One{ y_i-\yh_i>z_i>0}\\
&\geq\sum_i (y_i - \yh_i -z_i)_+ \cdot \One{z_i>0}\\
&\geq \sum_i (\phi_i^\top(x-\xh ) -2/(\sigma n)-z_i)_+ \cdot \One{z_i>0} \\
&\hspace{1in}{}- \sum_i |y_i - \phi_i^\top x| \cdot \One{z_i>0}\cdot \One{|y_i -\phi_i^\top x|> 2/(\sigma n)}\\
&\geq \sum_i (\phi_i^\top(x-\xh ) -z_i)_+ \cdot \One{z_i>0} - \frac{2}{\sigma n}\sum_i \One{z_i>0} - \frac{\sigma n}{2} \sum_i (y_i - \phi_i^\top x)^2 \cdot \One{z_i>0}
\end{align*}
and similarly,
\[\textnormal{Term 2} \geq \sum_i (-\phi_i^\top(x-\xh ) + z_i)_+ \cdot \One{z_i<0}  - \frac{2}{\sigma n}\sum_i \One{z_i<0} - \frac{\sigma n}{2}\sum_i (y_i - \phi_i^\top x)^2 \cdot \One{z_i<0}.\]
Therefore, defining
\begin{equation}\label{eqn:Hx}H(x) = \frac{1}{n}\sum_i\left[(\phi_i^\top x-z_i)_+ \cdot \One{z_i>0} +  (-\phi_i^\top x + z_i)_+ \cdot \One{z_i<0}  \right],\end{equation}
and simplifying, we have
\[\inner{y-\yh }{\zeta_y - \zeta_{\yh}} \geq H(x-\xh)- \frac{\sigma}{2}\norm{y-\Phi x}^2_2- \frac{2}{\sigma n}.\]
We will now use the following lemma (proved in Appendix~\ref{app:proof_lem:QR_Hx}):
\begin{lemma}\label{lem:QR_Hx}
Suppose that $n\geq 4$, and that $\phi_1,\dots,\phi_n\in\R^d$ and $z_1,\dots,z_n\in\R$ satisfy
the assumptions of Proposition~\ref{prop:quantile_rsc}.
Then, with probability at least $1-(2nd)^{-1}$,
\[H(x)\geq C^*_1 \cdot \min\{\norm{x}^2_2,\norm{x}_2\} - C^*_2 \cdot \sqrt{\frac{\log(nd)}{n}}\cdot \norm{x}_1 - C^*_3\cdot \frac{\sqrt{\log(nd)}}{n}\text{ for all $x\in\R^d$},\]
where $H(x)$ is defined as in~\eqref{eqn:Hx}, and where $C^*_1,C^*_2,C^*_3$ are positive and finite, and depend only on
the constants $a_\phi,b_\phi,B_\phi,c_z,t_z$ appearing in the assumptions of Proposition~\ref{prop:quantile_rsc}.
\end{lemma}
 Returning to our work above
we therefore see that, with probability at least $1-(2nd)^{-1}$,
\begin{multline*}\inner{y-\yh }{\zeta_y - \zeta_{\yh}} \geq C^*_1  \min\{\norm{x-\xh}^2_2,\norm{x-\xh}_2\}\\ - C^*_2  \sqrt{\frac{\log(nd)}{n}}\cdot \norm{x-\xh}_1- \frac{\sigma}{2}\norm{y-\Phi x}^2_2- \left( \frac{C^*_3 \sqrt{\log(nd)}}{n} + \frac{2}{\sigma n}\right)\end{multline*}
for all $x\in\dom(f)$, $y\in\R^n$, and $\zeta_y\in\partial g(y)$.

\paragraph{Combining the $x$ and $y$ terms}
Combining our bounds for the $x$ and $y$ terms, we have shown that, with probability at least $1-(nd)^{-1}$,
for all $x\in\dom(f)$, $y\in\R^n$, $\xi_x\in\partial f(x)$, and $\zeta_y\in\partial g(y)$,
\begin{multline*}
\biginner{\left(\begin{array}{c}x - \xh\\ y - \yh\end{array}\right)}{\left(\begin{array}{c} \xi_x - \xi_{\xh} \\ \zeta_y-\zeta_{\yh}\end{array}\right)} 
\geq \Bigg[\left(C_1 - 2B_\phi\right)\sqrt{\frac{\log(nd)}{n}}\norm{x-\xh}_1 
\\
 - \left(2C_2^{-1} + C_1\sqrt{C_0} \right) \cdot\min\{\norm{x-\xh}^2_2,\norm{x-\xh}_2\}
 -  \frac{C_2 C_\lambda^2}{4}\cdot\frac{s_*\log(nd)}{n}\Bigg]\\
+ \Bigg[C^*_1  \min\{\norm{x-\xh}^2_2,\norm{x-\xh}_2\}\\ - C^*_2  \sqrt{\frac{\log(nd)}{n}}\cdot \norm{x-\xh}_1- \frac{\sigma}{2}\norm{y-\Phi x}^2_2- \left( \frac{C^*_3 \sqrt{\log(nd)}}{n} + \frac{2}{\sigma n}\right)\Bigg].
\end{multline*}
We can simplify this to
\begin{multline*}
\biginner{\left(\begin{array}{c}x - \xh\\ y - \yh\end{array}\right)}{\left(\begin{array}{c} \xi_x - \xi_{\xh} \\ \zeta_y-\zeta_{\yh}\end{array}\right)} 
\geq \left(C_1 - 2B_\phi - C^*_2\right)\sqrt{\frac{\log(nd)}{n}}\norm{x-\xh}_1 
\\
 +(C^*_1 -  2C_2^{-1} - C_1\sqrt{C_0}) \min\{\norm{x-\xh}^2_2,\norm{x-\xh}_2\} \\
 - \frac{\sigma}{2}\norm{y-\Phi x}^2_2- \left( \frac{C^*_3 \sqrt{\log(nd)}}{n} + \frac{2}{\sigma n}+ \frac{C_2 C_\lambda^2}{4}\cdot \frac{s_*\log(nd)}{n}\right).
\end{multline*}
Choosing $C_1 = 2B_\phi+C^*_2$, , $C_2 = 8/C^*_1$, and $C_4 =C^*_3+ 2+ \frac{C_2 C_\lambda^2}{4}$,
and choosing $C_0$ to satisfy $C_0 \leq (C^*_1/4C_1)^2$, this simplifies to
\begin{multline*}
\biginner{\left(\begin{array}{c}x - \xh\\ y - \yh\end{array}\right)}{\left(\begin{array}{c} \xi_x - \xi_{\xh} \\ \zeta_y-\zeta_{\yh}\end{array}\right)} 
\geq \frac{C^*_1}{2} \min\{\norm{x-\xh}^2_2,\norm{x-\xh}_2\} \\
 - \frac{\sigma}{2}\norm{y-\Phi x}^2_2-  C_4\max\{1,\sigma^{-1}\}\cdot \frac{s_*\log(nd)}{n}.
\end{multline*}
To complete our proof that~\eqref{eqn:quantile_rsc_to_show} holds, we only need to verify that $C_1^*/2\geq C_3$.
Recall that we have defined this constant as $C_3 = 2\sqrt{C_0}(C_1 + 2B_\phi)$.
Therefore, by taking $C_0 =(C_1^*/4(C_1+2 B_\phi))^2$, all the necessary bounds are verified and we have completed the proof.

\subsubsection{Proof of Lemma~\ref{lem:QR_Hx}}\label{app:proof_lem:QR_Hx}

For any fixed $x$, define
\[\tilde{H}(x) = \EE{H(x)} = \EE{(\phi^\top x-z)_+ \cdot \One{z>0} +  (-\phi^\top x + z)_+ \cdot \One{z<0}},\]
where the expectation is taken with respect to $\phi\sim\mathcal{D}_\phi$ and $z\sim h_z$, with $\phi\independent z$.
We can calculate
\begin{align*}
&\EEst{(\phi^\top x-z)_+ \cdot \One{z>0}}{\phi}
=\int_{t=0}^{(\phi^\top x)_+} \big[(\phi^\top x)_+ - t\big] h_z(t)\;\mathsf{d}t\\
&\hspace{1in}\geq\int_{t=0}^{\min\{t_z,(\phi^\top x)_+\}} \big[(\phi^\top x)_+ - t\big] c_z\;\mathsf{d}t\\
&\hspace{1in}=\left[c_z\min\{(\phi^\top x)_+^2,t_z(\phi^\top x)_+\} - \frac{c_z}{2} \min\{(\phi^\top x)_+,t_z\}^2\right]\\
&\hspace{1in}\geq \frac{c_z}{2}\min\{(\phi^\top x)_+^2,t_z(\phi^\top x)_+\}.
\end{align*}
Similarly,
\[\EEst{(-\phi^\top x+z)_+ \cdot \One{z<0}}{\phi}\geq \frac{c_z}{2}\min\{(-\phi^\top x)_+^2,t_z(-\phi^\top x)_+\}.\]
Therefore, for a fixed $x$,
\begin{multline*}\tilde{H}(x) =\EE{(\phi^\top x-z)_+ \cdot \One{z>0} +  (-\phi^\top x + z)_+ \cdot \One{z<0}}\\
=\EE{\EEst{(\phi^\top x-z)_+ \cdot \One{z>0} +  (-\phi^\top x + z)_+ \cdot \One{z<0}}{\phi}}\\
\geq \EE{ \frac{c_z}{2}\min\{(\phi^\top x)^2,t_z|\phi^\top x|\}}.
\end{multline*}
Next, for any unit vector $u$, we can calculate	
\begin{multline*}
\EE{|\phi^\top u|^2\cdot \One{|\phi^\top u|\leq \frac{2b_\phi}{a_\phi}}}
=\EE{|\phi^\top u|^2} - \EE{|\phi^\top u|^2\cdot \One{|\phi^\top u|> \frac{2b_\phi}{a_\phi}}}\\
\geq\EE{|\phi^\top u|^2} - \frac{a_\phi}{2b_\phi}\EE{|\phi^\top u|^3}
\geq \frac{a_\phi}{2},
\end{multline*}
by our assumptions on $\mathcal{D}_\phi$.
For any $x\neq 0$, writing $u=\frac{x }{\norm{x }_2}$,
\begin{align*}
&\min\big\{(\phi^\top x)^2,t_z|\phi^\top x|\big\}\\
&\hspace{.5in}=\min\big\{\norm{ x}^2_2\cdot (\phi^\top u)^2,t_z \norm{ x}_2 \cdot |\phi^\top u|\big\}\\
&\hspace{.5in}\geq \min\left\{\norm{ x}^2_2,\frac{ t_za_\phi}{2b_\phi}\norm{ x}_2\right\}\cdot \min\left\{ (\phi^\top u)^2,\frac{2b_\phi}{a_\phi} |\phi^\top u|\right\}\\
&\hspace{.5in}\geq \min\left\{\norm{ x}^2_2,\frac{ t_za_\phi}{2b_\phi}\norm{ x}_2\right\}\cdot (\phi^\top u)^2 \cdot \One{|\phi^\top u|\leq \frac{2b_\phi}{a_\phi}},
\end{align*}
and therefore for all $x$,
\begin{multline*}
\EE{\min\big\{(\phi^\top x)^2,t_z|\phi^\top x|\big\}}\\
\geq \EE{\min\left\{\norm{ x}^2_2,\frac{a_\phi t_z}{2b_\phi}\norm{ x}_2\right\}\cdot (\phi^\top u)^2 \cdot \One{|\phi^\top u|\leq \frac{2b_\phi}{a_\phi}}}\\
\geq\frac{a_\phi}{2}\min\left\{\norm{ x }^2_2, \frac{t_za_\phi}{2b_\phi}\norm{ x }_2\right\}.\end{multline*}
Combining this with the work above,
\[\tilde{H}(x) \geq \frac{a_\phi c_z}{4}\min\left\{\norm{ x }^2_2, \frac{t_za_\phi}{2b_\phi}\norm{ x }_2\right\}\]
for all $x\in\R^d$.

Next, we will use a peeling argument to bound $\left|H(x) - \tilde{H}(x)\right|$.
 First, fixing any $B>0$,
\begin{multline*}
\EE{\sup_{x : \norm{x}_1\leq B} \left|H(x) - \tilde{H}(x)\right|}
\\\leq 2\EE{\sup_{x : \norm{x}_1\leq B} \left|\frac{1}{n}\sum_i \xi_i \left[(\phi_i^\top x-z_i)_+ \cdot \One{z_i>0} +  (-\phi_i^\top x + z_i)_+ \cdot \One{z_i<0}  \right]\right|}\end{multline*}
 by symmetrization \citep[Theorem 2.1]{koltchinskii2011oracle}, where $\xi_1,\dots,\xi_n\iidsim\textnormal{Unif}\{\pm 1\}$.
 Next, fixing the $z_i$'s, define $\varphi_i(t) = (t-z_i)_+\cdot\One{z_i>0} + (-t+z_i)_+\cdot\One{z_i<0}$. Then $\varphi_i$ is $1$-Lipschitz for all $i$, and so
 \begin{multline*}\EEst{\sup_{x : \norm{x}_1\leq B} \left|\frac{1}{n}\sum_i \xi_i \left[(\phi_i^\top x-z_i)_+ \cdot \One{z_i>0} +  (-\phi_i^\top x + z_i)_+ \cdot \One{z_i<0}  \right]\right|}{z_{1:n}}\\ \leq 2\EEst{\sup_{x : \norm{x}_1\leq B} \left|\frac{1}{n}\sum_i \xi_i \cdot \phi_i^\top x\right|}{z_{1:n}}
= 2\EE{\sup_{x : \norm{x}_1\leq B} \left|\frac{1}{n}\sum_i \xi_i \cdot \phi_i^\top x\right|}\end{multline*}
 by the Rademacher comparison inequality \citep[Theorem 2.2]{koltchinskii2011oracle}. Finally,
\[ 
\EE{\sup_{x : \norm{x}_1\leq B} \left|\frac{1}{n}\sum_i \xi_i \cdot \phi_i^\top x\right|}
 \leq \EE{\sup_{x : \norm{x}_1\leq B} \norm{\frac{1}{n}\Phi^\top \xi}_{\infty}\norm{ x }_1}
 = B\EE{\norm{\frac{1}{n}\Phi^\top \xi}_{\infty}} .\]
 And, we know that $\norm{\frac{1}{n}\Phi^\top \xi}_{\infty}\leq B_\phi$ deterministically (since $\norm{\xi}_{\infty}\leq 1$ and $\norm{\Phi}_{\infty}\leq B_\phi$),
 and so applying~\eqref{eqn:QR_Hoeffding}, we have
 \[\EE{\norm{\frac{1}{n}\Phi^\top \xi}_{\infty}} \leq 2B_\phi \sqrt{\frac{\log(nd)}{n}} + \frac{B_\phi}{2nd} \leq 3B_\phi \sqrt{\frac{\log(nd)}{n}},\]
where the last step holds since $\frac{1}{2nd}\leq  \sqrt{\frac{\log(nd)}{n}}$ for all $n\geq 4,d\geq 1$.
So, we have
\[ \EE{\sup_{x : \norm{x}_1\leq B} \left|\frac{1}{n}\sum_i \xi_i \cdot \phi_i^\top x\right|}
 \leq 3BB_\phi \sqrt{\frac{\log(nd)}{n}}.\]
 Combining our work so far, we have shown that
\[\EE{\sup_{x : \norm{x}_1\leq B} \left|H(x) - \tilde{H}(x)\right|} \leq 
4\EE{\sup_{x : \norm{x}_1\leq B} \left|\frac{1}{n}\sum_i \xi_i \cdot \phi_i^\top x\right|}\\
\leq 12BB_\phi\sqrt{\frac{\log(nd)}{n}}.\]

Next, if we alter one data point $\phi_i,z_i$, we can see that the value of $H(x)$ changes by at most $\frac{1}{n} B_\phi\norm{x}_1$. 
Therefore, by McDiarmid's inequality, for any $t>0$,
\[\PP{\sup_{x : \norm{x}_1\leq B} \left|H(x) - \tilde{H}(x)\right| > \EE{\sup_{x : \norm{x}_1\leq B} \left|H(x) - \tilde{H}(x)\right|} + t}\leq \exp\left\{ - \frac{2t^2}{n^{-1}B^2B_\phi^2}\right\}.\]
Setting $t = BB_\phi \sqrt{\log (nd)/n}$ and plugging in our calculation for the expected value,
\[\PP{\sup_{x : \norm{x}_1\leq B} \left|H(x) - \tilde{H}(x)\right| > 13BB_\phi\sqrt{\frac{\log (nd)}{n}}}\leq (nd)^{-2}.\]

 Therefore, applying this result with $B=\sqrt{d}, 2^{-1}\sqrt{d},\dots,2^{-(K-1)}\sqrt{d}$ for $K=1+ \lceil\frac{1}{2} \log_2(nd)\rceil$ (i.e., a peeling
 argument, with $K$ chosen so that the smallest value of $B$ is $\leq n^{-1/2}$),  we see that $K\leq nd/2$ (this holds for any $n\geq 4,d\geq 1$),
 and so with probability at least $1-(2nd)^{-1}$,
\[\left|H(x) - \tilde{H}(x)\right| \leq \max\{\norm{ x }_1 , n^{-1/2}\} \cdot  26B_\phi \sqrt{\frac{\log(nd)}{n}}\]
for all $x$ with $ \norm{ x }_1 \leq \sqrt{d}$---and therefore, for all $x$ with $\norm{ x }_2\leq 1$.
Combining everything so far, we have shown that with probability at least $1-(2nd)^{-1}$,
\begin{equation}\label{eqn:Hx_bound_unitball}H(x)\geq \frac{a_\phi c_z}{4}\min\left\{1, \frac{t_za_\phi}{2b_\phi}\right\}\cdot\norm{x}^2_2 - \max\{\norm{ x }_1 , n^{-1/2}\} \cdot  26B_\phi \sqrt{\frac{\log(nd)}{n}}\end{equation}
for all $x\in\R^d$ with $\norm{ x}_2\leq 1$.

Now we consider $x$ with $\norm{ x }_2\geq 1$. Let $x' =  \frac{ x}{\norm{ x}_2}$.
Since $H(x)$ is convex, and $1=\norm{x'}_2\leq\norm{x'}_1$, if the bound~\eqref{eqn:Hx_bound_unitball} holds (at $x=x'$) then we have
\[
\frac{1}{\norm{ x}_2} H(x) + \left(1-\frac{1}{\norm{ x}_2}\right)H(0) \geq H(x') \geq \frac{a_\phi c_z}{4}\min\left\{1, \frac{t_za_\phi}{2b_\phi}\right\} -  \norm{x' }_1 \cdot  26B_\phi \sqrt{\frac{\log(nd)}{n}}.\]
Clearly $H(0)=0$, and since $\norm{ x}_1 = \norm{x'}_1\cdot\norm{ x}_2$, we can simplify this to
\[H(x)\geq \frac{a_\phi c_z}{4}\min\left\{1, \frac{t_za_\phi}{2b_\phi}\right\}\norm{ x }_2 -  \norm{ x }_1 \cdot  26B_\phi \sqrt{\frac{\log(nd)}{n}}.\]
Combining both cases (i.e., $\norm{ x}_2\leq1$ or $\norm{ x}_2>1$), we have therefore proved that, with probability at least $1-(2nd)^{-1}$,
\[H(x)\geq \frac{a_\phi c_z}{4}\min\left\{1, \frac{t_za_\phi}{2b_\phi}\right\}\cdot\min\{\norm{x}^2_2,\norm{x}_2\} - \max\{\norm{ x }_1 , n^{-1/2}\} \cdot  26B_\phi \sqrt{\frac{\log(nd)}{n}},\]
for all $x\in\R^d$, which completes the proof.

\end{document}